\documentstyle[12pt]{article} 

\def\ds{\displaystyle}
\def\be{\begin{equation}} 
\def\ee{\end{equation}} 
\newtheorem{th}{Theorem}[section]  
\newtheorem{lm}[th]{Lemma}

\def\RR{{\hbox{I\kern-.2em\hbox{R}}}}
\def\PP{{\hbox{I\kern-.2em\hbox{P}}}}
\def\EE{{\hbox{I\kern-.2em\hbox{E}}}}
\def\ZZ{{\hbox{Z\kern-.4em\hbox{Z}}}}
\def\br{\mbox{\boldmath $r$}}
\def\bz{\mbox{\boldmath $z$}}
\def\one{{\bf 1}}

\def\R{{\bf R}}

\def\qed{\hfill\hbox{\rule{6pt}{6pt}}}
\def\essinf{\mathop{\rm ess \!~inf}_{(E,m)}}
\def\esssup{\mathop{\rm ess \!~sup}_{(E,m)}}
\def\cdotsi{\mathop{\check{\cdots}}^{i}}
\def\cdotsij{\mathop{\check{\cdots}}^{i,j}}
\def\lg{\langle}
\def\rg{\rangle}

\def\cA{{\cal A}}

\def\cD{{\cal D}}
\def\cE{{\cal E}}
\def\cF{{\cal F}}

\def\cM{{\cal M}}

\def\cS{{\cal S}}
\def\cV{{\cal V}}

\def\cL{{\cal L}}

\begin{document} 
\setlength{\topmargin}{-4mm}
%
%
\begin{center}
{\LARGE \bf 
Ergodic properties for 
$\alpha$-CIR models 
and a class of 
generalized Fleming-Viot processes 
}
\end{center} 

%
%


\begin{center}
Kenji Handa 
\end{center} 

\begin{center}
{ 
Department of Mathematics \\
Saga University \\ 
Saga 840-8502 \\ 
Japan
} \\  
e-mail: 
handa@ms.saga-u.ac.jp \\ 
FAX: +81-952-28-8501
\end{center} 
\begin{center}
{\em Dedicated to Professor Ken-iti Sato 
on the occasion of his 80th birthday}
\end{center}

%
%

\begin{center} 
\begin{minipage}[t]{12cm} 
\small 
We discuss a Markov jump process 
regarded as a variant of 
the CIR (Cox-Ingersoll-Ross) model 
and its infinite-dimensional extension. 
These models belong to a class of 
measure-valued branching processes with immigration, 
whose jump mechanisms are 
governed by certain stable laws. 
The main result gives a lower spectral gap estimate 
for the generator. 
As an application, a certain ergodic property is shown 
for the generalized Fleming-Viot process 
obtained as the time-changed ratio process. 
\end{minipage}
\end{center}

%
%

\section{Introduction} 
\setcounter{equation}{0} 
The study of ergodic behaviors of a Markov process is 
of quite interest for various reasons. 
For instance, it is typical that 
the analysis of such behaviors depends 
heavily on the mathematical structure of the model, 
so that resulting properties are 
expected to yield deep understanding for it. 
In this paper, we discuss two specific classes of 
measure-valued Markov jump processes. 
The one consists of what we will call 
measure-valued $\alpha$-CIR models, 
each of which is thought of 
as an infinite-dimensional extension 
for a jump-type version of the CIR model, 
and the other generalizes naturally 
a class of Fleming-Viot processes 
with parent-independent mutation. 
As for the measure-valued $\alpha$-CIR model, 
identification of a stationary distribution 
is easy thanks to its nice structure 
as a measure-valued branching process 
with immigration (henceforth MBI-process). 
For the latter class of models, 
stationary distributions 
are identified recently in \cite{H13}. 
A key idea there is 
to exploit a special relationship with 
measure-valued $\alpha$-CIR models, 
which enabled us to give an expression 
for stationary distributions of 
our generalized Fleming-Viot processes 
in terms of those of 
the measure-valued $\alpha$-CIR models. 
It should be mentioned that 
such links have been discussed in another context 
in \cite{BBC} and \cite{FH}. 
Our attempt here is to rely still on that relationship 
to explore ergodic properties for 
the generalized Fleming-Viot process.

It is worth illustrating by taking up 
a one-dimensional model which is regarded 
as a `prototype' of the above mentioned MBI-process. 
Consider the well-known 
CIR model governed by generator 
\be 
L_{1}=z\frac{d^2}{dz^2}+(-bz+c)\frac{d}{dz}, 
\quad z\in\R_+:=[0,\infty), 
                                      \label{1.1} 
\ee 
where $b\in\R$ and $c>0$ are constants. 
Rather than its importance 
in the context of mathematical finance, 
we emphasize that this model belongs to 
the class of continuous state branching processes 
with immigration (CBI-processes in short). 
(See \cite{KW} for fundamental results 
regarding this class.) 
Let $0<\alpha<1$ be arbitrary. 
As a natural non-local version of (\ref{1.1}) 
within the class of generators of 
conservative  CBI-processes 
(cf. Theorem 1.2 in \cite{KW}), 
we will be concerned with 
\begin{eqnarray}
{L}_{\alpha}F(z)
&=& 
\frac{\alpha+1}{\Gamma(1-\alpha)}
z\int_0^{\infty}\left[F(z+y)-F(z)-yF'(z)\right]
\frac{dy}{y^{\alpha+2}} 
\nonumber  \\ 
& & 
-\frac{b}{\alpha}zF'(z)+c\frac{\alpha}{\Gamma(1-\alpha)}
\int_0^{\infty}\left[F(z+y)-F(z)\right]
\frac{dy}{y^{\alpha+1}},   
                                             \label{1.2} 
\end{eqnarray} 
where $\Gamma(\cdot)$ is the gamma function. 
The operator $L_{\alpha}$ with $b=0$ is found 
in Example 1.1 of \cite{KW}. 
Observing that, as $\alpha\uparrow 1$, 
${L}_{\alpha}F(z)\to {L}_{1}F(z)$ 
for any  $z>0$ and `nice' functions $F$ on $\R_+$, 
we call a Markov process associated with $a{L}_{\alpha}$ 
for some constant $a>0$ an $\alpha$-CIR model.

Although this class of models would be 
of interest in its own right especially 
in the mathematical finance context, 
our main motivation to study it  
is the analysis of ergodicity for 
a jump-type version of 
a Wright-Fisher diffusion model with mutation, 
which is obtained through normalization and 
random time-change from two independent processes 
with generators of the form (\ref{1.2}), 
say $L_{\alpha}'$ and $L_{\alpha}''$, 
with common $\alpha$ and $b$.   
On the level of generators such a link can be 
reformulated as the identity 
\begin{eqnarray} 
\lefteqn{
\left(L_{\alpha}'F(\cdot,z_2)\right)(z_1)
+ \left(L_{\alpha}''F(z_1,\cdot)\right)(z_2)} 
\nonumber  \\ 
& = & C(z_1+z_2)^{-\alpha}
A_{\alpha}
G\left(\frac{z_1}{z_1+z_2}\right), 
\quad z_1,z_2>0, 
                                      \label{1.3} 
\end{eqnarray} 
where $G$ is any smooth function on $[0,1]$, 
$F$ is defined by $F(z_1,z_2)=G(z_1/(z_1+z_2))$,  
$C$ is a positive constant independent of $G$ 
and $A_{\alpha}$ is 
the generator of a jump-type version of 
the Wright-Fisher diffusion model. 
(See (1.3) in \cite{H13} for a concrete 
expression for $A_{\alpha}$ or 
(\ref{5.1}) below for its generalization.) 
A significant consequence of (\ref{1.3}) 
is that Dirichlet form associated with $A_{\alpha}$ is, 
up to some multiplicative constant, 
a restriction of Dirichlet form 
associated with the two independent 
$\alpha$-CIR models. Therefore, 
ergodic properties of $\alpha$-CIR models 
would be expected to help us obtain the same kind of 
results for the process associated with $A_{\alpha}$.

Such an idea can extend naturally to 
the measure-valued $\alpha$-CIR model, which is  
regarded roughly as `continuum direct sum' 
of $\alpha$-CIR models with 
coefficients depending on a spatial parameter. 
More precisely, it is a MBI-process on a type space $E$, say, 
with zero mutation, branching mechanism 
\begin{eqnarray} 
E\times \R_+ \ni (r,\lambda) 
& \mapsto & 
a(r)\frac{\alpha+1}{\Gamma(1-\alpha)}
\int_0^{\infty}\left[e^{-\lambda y}-1+\lambda y\right]
\frac{dy}{y^{\alpha+2}} 
-\frac{b(r)}{\alpha}\lambda    \nonumber \\ 
& & 
= -\frac{1}{\alpha}(a(r)\lambda^{\alpha+1}+b(r)\lambda) 
                                          \label{1.4} 
\end{eqnarray} 
and nonlocal immigration mechanism 
\be 
f(\cdot) \mapsto 
\frac{\alpha}{\Gamma(1-\alpha)}\int_E m(dr)
\int_0^{\infty}\left[1-e^{-f(r) y}\right]
\frac{dy}{y^{\alpha+2}} 
= \int_E m(dr)f(r)^{\alpha},             \label{1.5} 
\ee 
where $a(r)>0, b(r)\in \R$ and 
$m$ is a finite non-null measure on $E$. 
An important feature of such processes is that 
the transition semigroups (for fixed $t$) 
and stationary distributions 
form a convolution semigroup with respect to $m$. 
(See (\ref{2.15}) and (\ref{2.17}) below.) 
Because of this structure 
studying ergodic properties of the extended model 
would be reduced to the one-dimensional case 
at least under the assumption of uniform bounds 
for the coefficients. In addition, 
as observed in \cite{H13}, the relation (\ref{1.3}) 
admits a generalization 
in the setting of measure-valued processes. 
(See also (\ref{5.2}) below.) 
For this reason the above mentioned extension 
of the $\alpha$-CIR model is considered 
to play an important role in studying 
the generalized Fleming-Viot process 
obtained as the time-changed ratio process.

The organization of this paper is as follows. 
In Section 2, we introduce the 
measure-valued $\alpha$-CIR model, 
and it is shown in Section 3 
that a lower spectral gap estimate for the generator 
can reduce to the one-dimensional case 
in a suitable sense. 
In Section 4, we prove such an estimate for $L_{\alpha}$, 
establishing exponential convergence to equilibrium 
for the measure-valued $\alpha$-CIR model. 
The latter result will be applied to a class of 
generalized Fleming-Viot processes in Section 5. 
\section{The measure-valued $\alpha$-CIR models}
\setcounter{equation}{0} 
To discuss in the setting of measure-valued processes, 
we need the following notation. 
Let $E$ be a compact metric space and 
$C(E)$ (resp. $B_+(E)$) 
the set of continuous (resp. nonnegative, bounded Borel) 
functions on $E$. Also, 
denote by $C_{++}(E)$ the set of functions in $C(E)$ 
which are uniformly positive. 
Define $\cM(E)$ to be the totality of finite Borel measures on $E$, 
and we equip $\cM(E)$ with the weak topology. 
Denote by $\cM(E)^{\circ}$ the set of non-null elements of $\cM(E)$. 
The set $\cM_1(E)$ of Borel probability measures on $E$ is 
regarded as a subspace of $\cM(E)$. 
We also use notation $\lg \eta, f\rg:=\int_E f(r)\eta(dr)$.  
For each $r\in E$, let $\delta_r$ denote 
the delta distribution at $r$. Given a probability measure $Q$, 
we write also $E^Q[\cdot]$ for the expectation with 
respect to $Q$.

Suppose that $0<\alpha<1$, $a\in C_{++}(E)$, $b\in C(E)$
and $m\in\cM(E)^{\circ}$ are given. 
As a natural generalization of the $\alpha$-CIR model 
generated by (\ref{1.2}), 
we shall discuss in this section 
the Markov process on $\cM(E)$ associated with 
\begin{eqnarray} 
\lefteqn{\cL_{\alpha}\Psi(\eta)
=\cL_{\alpha}^{(1)}\Psi(\eta)
+\cL_{\alpha}^{(2)}\Psi(\eta)
+\cL_{\alpha}^{(3)}\Psi(\eta)} \nonumber  \\ 
& := & 
\frac{\alpha+1}{\Gamma(1-\alpha)}\int_0^{\infty}
\frac{dz}{z^{2+\alpha}}\int_E\eta(dr)a(r)
\left[\Psi(\eta+z\delta_r)-\Psi(\eta)
-z\frac{\delta\Psi}{\delta\eta}(r)\right]
-\frac{1}{\alpha}\lg \eta,b\frac{\delta\Psi}{\delta\eta}\rg                                                                       \nonumber  \\ 
&  & 
+\frac{\alpha}{\Gamma(1-\alpha)}\int_0^{\infty}
\frac{dz}{z^{1+\alpha}}\int_Em(dr)
\left[\Psi(\eta+z\delta_r)-\Psi(\eta)\right], 
\qquad \eta\in\cM(E),  
                                              \label{2.1}  
\end{eqnarray} 
where $\frac{\delta\Psi}{\delta\eta}(r)
=\left.\frac{d}{d\epsilon}
\Psi(\eta+\epsilon\delta_r)\right|_{\epsilon=0}$. 
The operator $\cL_{\alpha}^{(3)}$ describes 
the mechanism of immigration. 
(See (9.25) in \cite{Li} for a general form 
of generators of MBI-processes. In our model, 
there is no `motion process', 
whose generator is thus considered to be $A\equiv 0$.) 
Set $\Psi_f(\eta)=e^{-\langle \eta,f \rangle}$ 
for $f\in B_+(E)$ and define $\cD$ to be 
the linear span of functions 
$\Psi_{f}$ with $f\in C_{++}(S)$. 
It is immediate to see from (\ref{1.4}) and (\ref{1.5}) 
that for any $f\in B_+(E)$ 
\be 
\cL_{\alpha}\Psi_f(\eta)
=
\Psi_f(\eta)\frac{1}{\alpha}
\langle \eta,af^{\alpha+1}+bf \rangle
-\Psi_f(\eta)\langle m,f^{\alpha} \rangle. 
                                              \label{2.2}  
\ee 
$\cL_{\alpha}$ is well-defined also on  
the class $\cF$ of functions $\Psi$ of the form 
\be 
\Psi(\eta)=\varphi(\lg\eta,f_1\rg,\ldots,\lg\eta,f_n\rg) 
                                            \label{2.3}  
\ee 
for some $\varphi\in C_0^2(\R_+^n)$, $f_i\in C_{++}(E)$ 
and a positive integer $n$. 
Our first result below not only verifies this  
but also gives bounds for each 
$\cL_{\alpha}^{(k)}\Psi$ $(k\in\{1,2,3\})$ 
for a more general class of functions $\Psi$. 
In what follows $\Vert\cdot\Vert_{\infty}$ 
denotes the sup norm. 
Let $\widetilde{\cF}$ be the totality of functions 
$\Psi$ of the form (\ref{2.3}) 
with $\varphi\in C^2(\R_+^n)$ and 
$\mbox{\boldmath $f$}:=(f_1,\ldots,f_n)\in C_{++}(E)^n$ 
satisfying the following conditions; 
there exist nonnegative constants 
$C^{(i)}_j (1\le i,j\le n)$, 
$C_{k}^{(ij)} (1\le i,j,k\le n)$ 
and $\epsilon>0$ such that 
for each $i,j\in\{1,\ldots,n\}$ 
\be 
|\varphi_{i}(x_1,\ldots,x_n)|
\le \sum_{k=1}^n\frac{C^{(i)}_k}{x_k+\epsilon}
\quad \mbox{for any} \ (x_1,\ldots,x_n)
\in\R_{\mbox{\boldmath $f$}}^n
                                              \label{2.4}  
\ee 
and 
\be 
|\varphi_{ij}(x_1,\ldots,x_n)|
\le \sum_{k=1}^n\frac{C_{k}^{(ij)}}{(x_k+\epsilon)^2} 
\quad \mbox{for any} \ (x_1,\ldots,x_n)
\in\R_{\mbox{\boldmath $f$}}^n, 
                                              \label{2.5}  
\ee 
where $\ds{\varphi_{i}=\frac{\partial \varphi}{\partial x_i}}$, 
$\ds{\varphi_{ij}
=\frac{\partial^2 \varphi}{\partial x_i \partial x_j}}$ 
and $\R_{\mbox{\boldmath $f$}}^n$ is defined to be 
\[ 
\left\{(x_1,\ldots,x_n)\in(0,\infty)^n:~
\frac{\inf_{x\in E}f_i(x)}{\Vert f_j\Vert_{\infty}}
\le \frac{x_i}{x_j} 
\le \frac{\Vert f_i\Vert_{\infty}}{\inf_{x\in E}f_j(x)}
\ (1\le i,j\le n) \right\}. 
\] 
Note that 
$\lg\eta,\mbox{\boldmath $f$}\rg 
:=(\lg\eta,f_1\rg,\ldots,\lg\eta,f_n\rg)\in 
\R_{\mbox{\boldmath $f$}}^n$ for any $\eta\in\cM(E)^{\circ}$. 
Intuitively, these conditions enable one to control 
the effect of long-range jumps governed by stable laws, 
and are inspired by the calculations in 
the proof of Proposition 3.4 in \cite{H13}. 

\medskip 

\noindent 
{\it Example.}~
It will turn out in Section 5 that 
an important example of functions 
in $\widetilde{\cF}\setminus \cF$ is 
\[ 
\Psi(\eta)=\lg\eta,f_1\rg\cdots\lg\eta,f_n
\rg(\lg\eta,f_{n+1}\rg+\epsilon)^{-n}, 
\] 
where $f_i\in C_{++}(E)$, $\epsilon>0$ 
and $n$ is a positive integer. 
This function corresponds to 
$\varphi(x_1,\ldots,x_{n+1})=x_1\cdots x_{n}
(x_{n+1}+\epsilon)^{-n}$, for which 
the following are verified to hold: 
\[ 
\varphi_{i}(x_1,\ldots,x_{n+1}) 
=\left\{
\begin{array}{ll} 
\ds{x_1\cdotsi x_{n}(x_{n+1}+\epsilon)^{-n}} 
& (i\in\{1,\ldots,n\}) \\ 
-nx_1\cdots x_{n}(x_{n+1}+\epsilon)^{-(n+1)} 
& (i=n+1) 
\end{array} 
\right.  
\] 
and 
\[ 
\varphi_{ij}(x_1,\ldots,x_{n+1})
=\left\{
\begin{array}{ll} 
0 & (i=j\in\{1,\ldots,n\}) \\ 
\ds{x_1\cdotsij x_{n}(x_{n+1}+\epsilon)^{-n}} 
& (i,j\in\{1,\ldots,n\}, i\ne j) \\ 
\ds{-nx_1\cdotsi x_{n}(x_{n+1}+\epsilon)^{-(n+1)}} 
& (i\in\{1,\ldots,n\}, j=n+1) \\ 
\ds{n(n+1)x_1\cdots x_{n}(x_{n+1}+\epsilon)^{-(n+2)}} 
& (i=j=n+1).   
\end{array} 
\right.  
\] 
Here, $\ds{\cdotsi}$ (resp. $\ds{\cdotsij}$) indicates 
deletion of the $i$th (resp. $i$th and $j$th) factor(s). 
These equalities are sufficient to show inequalities 
of the form (\ref{2.4}) and (\ref{2.5}). 
We can take in particular $C_k^{(i)}=0=C_k^{(ij)}$ 
for any $i,j\in\{1,\ldots,n+1\}$ and $k\in\{1,\ldots,n\}$. 
\begin{lm}
(i) It holds that $\cF\subset \widetilde{\cF}$. \\ 
(ii) Let $\Psi\in \widetilde{\cF}$ be expressed as (\ref{2.3}) 
with $\varphi$ satisfying (\ref{2.4}) and (\ref{2.5}) 
and $f_i\in C_{++}(E)$. 
Then for any $\eta\in \cM(E)$ 
\be 
\cL_{\alpha}^{(1)}\Psi(\eta) 
= 
\frac{1}{\alpha\Gamma(1-\alpha)}\int_E\eta(dr)a(r)
\int_0^{\infty}u^{-\alpha}du
\sum_{i,j=1}^nf_i(r)f_j(r)
\varphi_{ij}(\lg\eta+u\delta_r,\mbox{\boldmath $f$}\rg),  
                                              \label{2.6} 
\ee 
\be 
\cL_{\alpha}^{(2)}\Psi(\eta) 
=
-\frac{1}{\alpha}\sum_{i=1}^n\lg\eta,bf_i\rg
\varphi_{i}(\lg\eta,\mbox{\boldmath $f$}\rg)  
                                              \label{2.7} 
\ee 
and 
\be  
\cL_{\alpha}^{(3)}\Psi(\eta)
=
\frac{1}{\Gamma(1-\alpha)}\int_Em(dr)
\int_0^{\infty}w^{-\alpha}dw 
\sum_{i=1}^nf_i(r)
\varphi_{i}(\lg\eta+w\delta_r,\mbox{\boldmath $f$}\rg). 
                                              \label{2.8} 
\ee  
Also, we have the bounds 
\begin{eqnarray} 
|\cL_{\alpha}^{(1)}\Psi(\eta)|  
& \le & 
\Gamma(\alpha)\sum_{i,j,k=1}^nC_{k}^{(ij)}
\frac{\lg \eta,af_if_jf_k^{\alpha-1}\rg}
{(\lg \eta,f_k\rg+\epsilon)^{\alpha+1}}        \nonumber \\ 
& \le & 
\Gamma(\alpha)\sum_{i,j,k=1}^nC_{k}^{(ij)}
\frac{\Vert af_if_jf_k^{\alpha-1}\Vert_{\infty}}
{\inf_{x\in E}f_k(x)}
(\lg \eta,f_k\rg+\epsilon)^{-\alpha},  
                                              \label{2.9}  
\end{eqnarray} 
\be 
|\cL_{\alpha}^{(2)}\Psi(\eta)| 
\ \le \ 
\frac{1}{\alpha}\sum_{i,j=1}^nC^{(i)}_j
\frac{\lg\eta,|b|f_i\rg}{\lg\eta,f_j\rg+\epsilon}  
\ \le \ 
\frac{1}{\alpha}\sum_{i,j=1}^nC^{(i)}_j
\frac{\Vert bf_i\Vert_{\infty}}{\inf_{x\in E}f_j(x)}  
                                              \label{2.10}  
\ee 
and 
\be 
|\cL_{\alpha}^{(3)}\Psi(\eta)| 
\le 
\Gamma(\alpha)\sum_{i,j=1}^nC^{(i)}_j
\frac{\lg m,f_if_j^{\alpha-1}\rg}
{(\lg \eta,f_j\rg+\epsilon)^{\alpha}}.    
                                              \label{2.11}  
\ee 
In particular, $\cL_{\alpha}^{(1)}\Psi$, 
$\cL_{\alpha}^{(2)}\Psi$ and 
$\cL_{\alpha}^{(3)}\Psi$ are bounded. 
\end{lm}
{\it Proof.}~ 
(i) Let $\varphi\in C_0^2(\R_+^n)$ be given 
and take $R_1,\ldots,R_n>0$ large enough so that 
$\varphi(x_1,\ldots,x_n)=0$ 
whenever $\max\{x_1/R_1,\ldots,x_n/R_n\}>1$. 
Let $\epsilon>0$ be arbitrary. Then it is easy 
to see that for all $(x_1,\ldots,x_n)\in\R_+^n$ 
\[ 
|\varphi_{i}(x_1,\ldots,x_n)|
\le 
\frac{1}{n}\sum_{k=1}^n
\frac{R_k+\epsilon}{x_k+\epsilon}\Vert 
\varphi_{i}\Vert_{\infty}
\] 
and 
\[
|\varphi_{ij}(x_1,\ldots,x_n)|
\le 
\frac{1}{n}\sum_{k=1}^n
\frac{(R_k+\epsilon)^2}{(x_k+\epsilon)^2}\Vert 
\varphi_{ij}\Vert_{\infty}. 
\] 
In view of (\ref{2.4}) and (\ref{2.5}), 
what we have just seen suffice to imply that 
$\cF\subset \widetilde{\cF}$. \\ 
(ii) 
First, we consider $\cL_{\alpha}^{(2)}\Psi(\eta)$, 
assuming that $\eta\in\cM(E)^{\circ}$. 
(If $\eta$ is the null measure, (\ref{2.10}) is trivial.) 
Observe that 
\be  
\frac{\delta\Psi}{\delta\eta}(r)
=\sum_{i=1}^nf_i(r)
\varphi_{i}(\lg\eta,\mbox{\boldmath $f$}\rg), 
                                              \label{2.12}  
\ee 
from which (\ref{2.7}) follows. 
Also, (\ref{2.10}) is immediate from (\ref{2.4}). 

The next task is to prove the assertions 
for $\cL_{\alpha}^{(3)}\Psi(\eta)$. 
Since $\frac{d}{dz}\Psi(\eta+z\delta_r)=
\frac{\delta\Psi}{\delta(\eta+z\delta_r)}(r)$, 
we have by Fubini's theorem 
\begin{eqnarray} 
\int_0^{\infty}
\frac{dz}{z^{1+\alpha}}
\left[\Psi(\eta+z\delta_r)-\Psi(\eta)\right] 
& = & 
\int_0^{\infty}
\frac{dz}{z^{1+\alpha}}
\int_0^zdw\frac{\delta\Psi}{\delta(\eta+w\delta_r)}(r)
                                               \nonumber  \\ 
& = & 
\frac{1}{\alpha}\int_0^{\infty}
w^{-\alpha}dw\frac{\delta\Psi}{\delta(\eta+w\delta_r)}(r).    
                                              \label{2.13} 
\end{eqnarray} 
So (\ref{2.8}) is deduced from (\ref{2.12}). 
Noting that $\eta+w\delta_r\in\cM(E)^{\circ}$ for $w>0$,   
apply (\ref{2.4}) to get 
\begin{eqnarray*} 
|\cL_{\alpha}^{(3)}\Psi(\eta)| 
& \le & 
\frac{1}{\Gamma(1-\alpha)}\int_Em(dr)
\int_0^{\infty}w^{-\alpha}dw\sum_{i,j=1}^n
\frac{f_i(r)C^{(i)}_j}{\lg\eta,f_j\rg+wf_j(r)+\epsilon} 
                                                      \\ 
& = & 
\Gamma(\alpha)\sum_{i,j=1}^nC^{(i)}_j\lg m,f_if_j^{\alpha-1}\rg 
(\lg \eta,f_j\rg+\epsilon)^{-\alpha}, 
\end{eqnarray*} 
which proves (\ref{2.11}). In the above equality we have used 
\be  
\int_0^{\infty}w^{-\alpha}\frac{dw}{sw+t} 
=\Gamma(\alpha)\Gamma(1-\alpha)s^{\alpha-1}t^{-\alpha},  
\quad s,t>0. 
                                              \label{2.14}  
\ee 

It remains to prove (\ref{2.6}) and (\ref{2.9}). 
Similarly to (\ref{2.13}) 
\begin{eqnarray*} 
\cL_{\alpha}^{(1)}\Psi(\eta)  
& = & 
\frac{\alpha+1}{\Gamma(1-\alpha)}\int_E\eta(dr)a(r)
\int_0^{\infty}\frac{dz}{z^{2+\alpha}}\int_0^zdw
\left[\frac{\delta\Psi}{\delta(\eta+w\delta_r)}(r)
-\frac{\delta\Psi}{\delta\eta}(r)\right] 
                                                     \\ 
& = & 
\frac{1}{\Gamma(1-\alpha)}\int_E\eta(dr)a(r)
\int_0^{\infty}\frac{dw}{w^{1+\alpha}}
\left[\frac{\delta\Psi}{\delta(\eta+w\delta_r)}(r)
-\frac{\delta\Psi}{\delta\eta}(r)\right] 
\end{eqnarray*} 
and by  (\ref{2.12}) 
\[ 
\frac{\delta\Psi}{\delta(\eta+w\delta_r)}(r)
-\frac{\delta\Psi}{\delta\eta}(r) 
=\int_0^wdu\sum_{i,j=1}^nf_i(r)f_j(r)
\varphi_{ij}(\lg\eta+u\delta_r,\mbox{\boldmath $f$}\rg). 
\] 
Hence (\ref{2.6}) is derived by Fubini's theorem. 
(\ref{2.6}) and (\ref{2.5}) together yield 
\begin{eqnarray*} 
|\cL_{\alpha}^{(1)}\Psi(\eta)| 
& \le & 
\frac{1}{\alpha\Gamma(1-\alpha)}\int_E\eta(dr)a(r)
\int_0^{\infty}u^{-\alpha}du                           
\sum_{i,j,k=1}^n\frac{f_i(r)f_j(r)C_{k}^{(ij)}}
{(\lg\eta,f_k\rg+uf_k(r)+\epsilon)^2}                      \\ 
& = & 
\Gamma(\alpha)\sum_{i,j,k=1}^nC_{k}^{(ij)}
\frac{\lg \eta,af_if_jf_k^{\alpha-1}\rg}
{(\lg \eta,f_k\rg+\epsilon)^{\alpha+1}}.  
\end{eqnarray*} 
Here, the last equality is deduced from 
\[ 
\int_0^{\infty}u^{-\alpha}\frac{du}{(su+t)^2} 
=\alpha\Gamma(\alpha)\Gamma(1-\alpha)
s^{\alpha-1}t^{-(\alpha+1)}, 
\quad s,t>0, 
\] 
which is verified by differentiating (\ref{2.14}) 
in $t$. 
\qed 

\medskip 

\noindent 
Following \cite{Stannat03}, 
we consider the operator $(\cL_{\alpha},\cF)$ as 
an operator on $C_{\infty}(\cM(E))$, the set of 
continuous functions on $\cM(E)$ vanishing at infinity. 
In the theorem below we collect basic properties 
of $\cL_{\alpha}$ and the associated transition semigroup. 
\begin{th}
(i) $(\cL_{\alpha},\cF)$ is closable in $C_{\infty}(\cM(E))$  
and the closure 
$(\overline{\cL_{\alpha}},D(\overline{\cL_{\alpha}}))$ 
generates a $C_0$-semigroup $(T(t))_{t\ge 0}$. 
Moreover, $\cD$ is a core for $\overline{\cL_{\alpha}}$, 
and for each $f\in B_+(E)$ and $\eta\in \cM(E)$  
\be 
T(t)\Psi_f(\eta) 
=\exp\left[-\lg \eta, V_tf\rg
-\int_0^t\lg m, (V_sf)^{\alpha}\rg ds\right], 
\qquad t\ge 0, 
                                               \label{2.15}  
\ee 
where 
\be   
V_tf(r)
=\frac{e^{-b(r)t/\alpha}f(r)}
{\left[1+a(r)f(r)^{\alpha}\int_0^te^{-b(r)s}ds 
\right]^{1/\alpha}}.   
                                               \label{2.16}  
\ee 
(ii) If $b\in C_{++}(E)$, then 
Markov process with transition semigroup 
$(T(t))_{t\ge 0}$ is ergodic in the sense that 
for every initial state $\eta\in \cM(E)$, 
the law of the process at time $t$ 
converges to a unique stationary distribution, 
say ${Q}_{\alpha}$, as $t\to\infty$. 
Moreover, the Laplace functional of ${Q}_{\alpha}$ is given by 
\be 
\int_{\cM(E)^{\circ}}{Q}_{\alpha}(d\eta)
\Psi_f(\eta) 
=\exp\left[-\langle m,a^{-1}\log(1+ab^{-1}f^{\alpha}) 
\rangle\right], 
\qquad f\in B_+(E). 
                                               \label{2.17}  
\ee 
\end{th}
{\it Proof.}~ 
(i) If $m$ were the null measure, 
the assertions except (\ref{2.16}) follow 
from more general Theorem 1.1 in \cite{Stannat03}, 
and also (\ref{2.16}) is deduced from the proof of it. 
Indeed, $V_tf(r)$ was given there implicitly by 
\be 
\frac{\partial}{\partial t}V_tf(r)
=-\frac{a(r)}{\alpha}V_tf(r)^{1+\alpha}
-\frac{b(r)}{\alpha}V_tf(r), 
\qquad V_0f(r)=f(r), 
                                           \label{2.18}  
\ee 
from which (\ref{2.16}) is obtained. 
(See Example 3.1 in \cite{Li}.) 
 
Based on these facts, 
the proof of the assertions for $m\in\cM(E)^{\circ}$ 
can be done by modifying suitably 
the proof of Corollary 1.3 in \cite{Stannat03}, 
which deals with the immigration mechanism 
described by the operator 
$\Psi \mapsto 
\lg m,\frac{\delta\Psi}{\delta\eta}\rg$. 
A (possibly unique) non-trivial modification would be  
the step to construct, for each $\eta\in\cM(E)$ 
and $t\ge0$, 
$q_t(\eta,\cdot)\in\cM_1(\cM(E))$ 
with Laplace transform given 
by the right side of (\ref{2.15}). 
By the observation made in the last paragraph, 
we have $p_t(\eta,\cdot)\in\cM_1(\cM(E))$ such that 
\[ 
\int_{\cM(E)}p_t(\eta,d\eta')\Psi_f(\eta')
=\exp\left[-\lg \eta, V_tf\rg\right], 
\qquad f\in B_+(E).  
\] 
Additionally, for every $\eta\in \cM(E)$, 
let $s_{\alpha}(\eta,\cdot)$ be the law of an 
$\alpha$-stable random measure with parameter 
measure $\eta$, i.e., 
\[ 
\int_{\cM(E)}s_{\alpha}(\eta,d\eta')\Psi_f(\eta')
=\exp\left[-\lg \eta, f^{\alpha}\rg\right], 
\qquad f\in B_+(E) 
\] 
and define $p_{t,\alpha}(\eta,\cdot)\in\cM_1(\cM(E))$ 
to be the mixture 
\[ 
p_{t,\alpha}(\eta,\cdot)
=\int_{\cM(E)}s_{\alpha}(\eta,d\eta')p_{t}(\eta',\cdot).  
\] 
It then follows that 
\[ 
\int_{\cM(E)}p_{t,\alpha}(\eta,d\eta')\Psi_f(\eta')
=\exp\left[-\lg \eta, (V_tf)^{\alpha}\rg\right], 
\qquad f\in B_+(E). 
\] 
Therefore, for each $N=1,2,\ldots$,  
the convolution 
\[ 
q_t^{(N)}(\eta,\cdot):=p_{t}(\eta,\cdot)
{\LARGE \ast} \left(\mathop{\ast}_{k=1}^N 
p_{tk/N,\alpha}\left(\frac{t}{N}m,\cdot\right)\right)
\] 
has Laplace transform 
\[ 
\int_{\cM(E)}q_t^{(N)}(\eta,d\eta')\Psi_f(\eta')
=\exp\left[-\lg \eta, V_tf\rg
-\sum_{k=1}^N\frac{t}{N}\lg m, 
(V_{tk/N}f)^{\alpha}\rg \right], 
\] 
which converges to the right side of 
(\ref{2.15}) as $N\to\infty$. Thus, the weak limit of 
$q_t^{(N)}(\eta,\cdot)$ as $N\to\infty$ is identified with 
the desired probability measure $q_t(\eta,\cdot)$ 
on $\cM(E)$.  
Hence the semigroup $(T(t))_{t \ge 0}$ defined by 
\[ 
T(t)\Psi(\eta)=\int_{\cM(E)}q_t(\eta,d\eta')\Psi(\eta'),  
\quad \Psi\in B(\cM(E)) 
\] 
satisfies (\ref{2.15}). The identity 
$\left.\frac{d}{dt}T(t)\Psi_f\right|_{t=0}=
\cL_{\alpha}\Psi_f$ for $f\in C_{++}(E)$ 
is verified by combining (\ref{2.2}) with (\ref{2.18}). 
Once (\ref{2.15}) is in hand, the assertion 
that $\cD$ is a core for $\overline{\cL_{\alpha}}$ 
follows as a direct consequence of Lemma 2.2 in \cite{Watanabe}. 
\\ 
(ii) 
As $t\to\infty$ the right side of (\ref{2.15}) converges to 
\[ 
\exp\left[-\int_0^{\infty}\lg m, (V_tf)^{\alpha}\rg dt\right]
= \exp\left[-\lg m, a^{-1}
\log(1+ab^{-1}f^{\alpha})\rg\right] 
\] 
since by (\ref{2.18}) 
\[ 
\frac{d}{dt}\log\left[1+a(r)b(r)^{-1}(V_tf(r))^{\alpha}\right] 
=-a(r)(V_tf(r))^{\alpha}. 
\] 
This proves the required ergodicity and that 
the unique stationary distribution $Q_{\alpha}$ 
has the Laplace functional given by 
the right side of (\ref{2.17}). 
The fact that $Q_{\alpha}$ is supported 
on $\cM(E)^{\circ}$ follows by observing that 
the right side of (\ref{2.17}) with 
$f\equiv \beta>0$ tends to 0 as $\beta\to \infty$.  
\qed 

\bigskip

\noindent 
We call the Markov process on $\cM(E)$ 
associated with (\ref{2.1}) in the sense of Theorem 2.2 
the measure-valued $\alpha$-CIR model 
with triplet $(a,b,m)$. 
It is said to be ergodic if $b\in C_{++}(E)$. 
\\ 
{\it Remarks.}~ 
(i) A random measure with law ${Q}_{\alpha}$ 
in Theorem 2.2 (ii) 
is an infinite-dimensional analogue of the 
random variable with law sometimes referred to as 
a (non-symmetric) Linnik distribution, 
whose Laplace exponent is of the form 
$\lambda\mapsto c\log(1+d\lambda^{\alpha})$ 
for some $c,d>0$. 
Observe from (\ref{2.17}) that, 
as $\alpha\uparrow 1$, ${Q}_{\alpha}$ 
converges to $Q_1$, the law of  
a generalized gamma process such that 
\[ 
\int_{\cM(E)^{\circ}}Q_1(d\eta)\Psi_f(\eta)
=\exp\left[-\langle m,a^{-1}\log(1+ab^{-1}f) 
\rangle\right], 
\qquad f\in B_+(E). 
\]  
In addition, one can see that 
\[ 
\lim_{\alpha\uparrow 1}\cL_{\alpha}\Psi(\eta) 
=\lg \eta,a\frac{\delta^2\Psi}{\delta\eta^2}\rg 
-\lg \eta,b\frac{\delta\Psi}{\delta\eta}\rg
+\lg m,\frac{\delta\Psi}{\delta\eta}\rg
=:\cL_1\Psi(\eta)   
\] 
for `nice' functions $\Psi$, 
where $\frac{\delta^2\Psi}{\delta\eta^2}(r)
=\left.\frac{d^2}{d\epsilon^2}
\Psi(\eta+\epsilon\delta_r)\right|_{\epsilon=0}$. 
(For instance, this is immediate for 
$\Psi=\Psi_f$ from (\ref{2.2}).) 
$\cL_1$ is the generator of an MBI-process 
discussed in Section 4 of \cite{Stannat03'} 
and in Section 3 of \cite{Stannat03},   
where $Q_{1}$ was shown to be 
a reversible stationary distribution 
of the process associated with $\cL_1$. \\ 
(ii) 
In contrast, ${Q}_{\alpha}$ ($0<\alpha<1$) 
is not a reversible stationary distribution of 
the measure-valued $\alpha$-CIR model. 
See Theorem 2.3 in \cite{H12} for 
an assertion of this type regarding CBI-processes. 
Essentially the same proof works at least 
in the case of ergodic measure-valued $\alpha$-CIR models. 
Namely, one can show, by a proof by contradiction, 
that the formal symmetry 
$E^{{Q}_{\alpha}}
\left[(-\cL_{\alpha})\Psi_f\cdot\Psi_g\right]
=E^{{Q}_{\alpha}}
\left[(-\cL_{\alpha})\Psi_g\cdot\Psi_f\right]$ 
fails for some $f,g\in C_{++}(E)$. 
For this purpose, an expression for 
the Dirichlet form $E^{{Q}_{\alpha}}
\left[(-\cL_{\alpha})\Psi_f\cdot\Psi_g\right]$ 
given Remark after Lemma 3.1 below is helpful. 
\section{Associated Dirichlet forms}
\setcounter{equation}{0} 
From now on, we suppose additionally that $b\in C_{++}(E)$. 
Thus, only ergodic measure-valued 
$\alpha$-CIR models will be discussed.  
To study the speed of convergence to equilibrium 
in the $L^2$-sense, we consider the symmetric part of 
Dirichlet form associated with $\cL_{\alpha}$ in (\ref{2.1}). 
It is a bilinear form on $\cF\times\cF$ 
defined by $\widetilde{\cE}(\Psi,\Psi')
:=E^{Q_{\alpha}}\left[\Gamma(\Psi,\Psi')\right]$ 
with $\Gamma(\cdot,\ast)$ being 
the `carr\'e du champ': 
\begin{eqnarray*}
\Gamma(\Psi,\Psi')(\eta) 
& := & 
\frac{1}{2}\left[-\Psi(\eta)\cL_{\alpha}\Psi'(\eta)
-\Psi'(\eta)\cL_{\alpha}\Psi(\eta)
+\cL_{\alpha}(\Psi\Psi')(\eta)\right]        \\ 
& =& 
\frac{1}{2}\int_{0}^{\infty}n_B(dz)\int_E\eta(dr)a(r)
\left[\Psi(\eta+z\delta_r)-\Psi(\eta)\right]
\left[\Psi'(\eta+z\delta_r)-\Psi'(\eta)\right] \\ 
& & 
+\frac{1}{2}\int_{0}^{\infty}n_I(dz)\int_Em(dr)
\left[\Psi(\eta+z\delta_r)-\Psi(\eta)\right]
\left[\Psi'(\eta+z\delta_r)-\Psi'(\eta)\right], 
\end{eqnarray*} 
where $n_B(dz)=(\alpha+1)z^{-\alpha-2}dz/\Gamma(1-\alpha)$ 
and $n_I(dz)=\alpha z^{-\alpha-1}dz/\Gamma(1-\alpha)$  
govern the jump mechanisms associated with 
branching and immigration, respectively. 
The same argument as in the proof of Proposition 1.6 
in \cite{Stannat03} shows that 
$(\cL_{\alpha},\cF)$ is closable in $L^2(Q_{\alpha})$ 
and that the closure 
$(\overline{\cL_{\alpha}}^{(2)},
D(\overline{\cL_{\alpha}}^{(2)}))$ 
generates a $C_0$-semigroup 
$({T}^2(t))_{t\ge0}$ on $L^2(Q_{\alpha})$ 
which coincides with $(T(t))_{t\ge0}$ 
when restricted to $C_{\infty}(\cM(E))$. 
We set $\cE(\Psi)=E^{Q_{\alpha}}
\left[(-\overline{\cL_{\alpha}}^{(2)})\Psi\cdot\Psi\right]$ 
for any $\Psi \in D(\overline{\cL_{\alpha}}^{(2)})$, 
remarking that $\cE(\Psi)=\widetilde{\cE}(\Psi,\Psi)$ 
for $\Psi\in\cF$.

Let ${\rm var}(\Psi)$ stand for the 
variance of $\Psi\in L^2(Q_{\alpha})$ 
with respect to $Q_{\alpha}$, namely, 
\[ 
{\rm var}(\Psi)
=E^{Q_{\alpha}}
\left[\left(\Psi-E^{Q_{\alpha}}[\Psi]\right)^2\right]. 
\] It is known that the largest $\kappa\ge 0$ such that 
\[ 
{\rm var}(T^2(t)\Psi)\le e^{-\kappa t}{\rm var}(\Psi)  
\quad \mbox{for all} \ 
\Psi\in L^2(Q_{\alpha}) \ \mbox{and} \ t> 0
\] 
is identified with 
\begin{eqnarray*} 
{\rm gap}(\overline{\cL_{\alpha}}^{(2)}) 
&:= & 
\inf\left\{\cE(\Psi):~ {\rm var}(\Psi)=1, 
\ \Psi\in D(\overline{\cL_{\alpha}}^{(2)})\right\}       \\ 
& = & 
\sup\left\{\kappa\ge 0:~  
\kappa\cdot{\rm var}(\Psi) \le \cE(\Psi) \ 
\mbox{for all} \ 
\Psi\in D(\overline{\cL_{\alpha}}^{(2)}) \right\}. 
\end{eqnarray*} 
We refer the reader to e.g. Theorem 2.3 in \cite{Liggett89} 
for the proof of this fact in a general setting. 
Besides, an estimate of the form 
${\rm gap}(\overline{\cL_{\alpha}}^{(2)})\ge \kappa$ 
implies that $\overline{\cL_{\alpha}}^{(2)}$ has 
a spectral gap below 0 of size 
larger than or equal to $\kappa$. 
(See Remark 1.13 in \cite{Stannat03}.) 
In calculating Dirichlet form 
and the variance functional with respect to $Q_{\alpha}$, 
we will make an essential use of the following expression 
for the `log-Laplace functional' in (\ref{2.17}): 
\begin{eqnarray} 
\psi(f)
& := & \langle m,a^{-1}\log(1+ab^{-1}f^{\alpha}) \rangle 
                                       \nonumber \\ 
& = & \int_Em_a(dr)\int_0^{\infty}
\Lambda(dz)\left(1-e^{-f^*(r)z}\right),     
                                               \label{3.1}  
\end{eqnarray} 
where $m_a(dr)=a(r)^{-1}m(dr)$, 
$\Lambda$ is the L\'evy measure of 
the infinite divisible distribution on $(0,\infty)$ 
with Laplace exponent 
$\lambda\mapsto \log(1+\lambda^{\alpha})$ 
and $f^*=(a/b)^{1/\alpha}f$. 
In what follows the domain of integration is 
understood to be $(0,\infty)$ when suppressed. 
Define nonnegative functions $K_B$ and $K_I$ 
on $\R_+^2$ by 
\be 
K_B(s ,t) 
:=  
\int n_B(dy)(1-e^{-sy})(1-e^{-ty}) 
=   
{\alpha}^{-1}\left[(s+t)^{\alpha+1}
-s^{\alpha+1}-t^{\alpha+1}\right]  
                                            \label{3.2} 
\ee 
and 
\be 
K_I(s,t) 
:= 
\int n_I(dy)(1-e^{-sy})(1-e^{-ty})
= 
s^{\alpha}+t^{\alpha}-(s+t)^{\alpha}, 
                                            \label{3.3} 
\ee 
respectively. The above identities are 
verified easily by differentiating in $s$ and $t$. 
\begin{lm}
For any $f,g\in B_+(E)$ 
\begin{eqnarray*} 
\widetilde{\cE}(\Psi_f,\Psi_g)
& = & 
\frac{1}{2}e^{-\psi(f+g)}\int_E m(dr)
\frac{a(r)^{1/\alpha}}{b(r)^{1/\alpha}}
\int n_B(dy)(1-e^{-{f}(r)y})(1-e^{-{g}(r)y})      \\ 
&   & 
\int\Lambda(dz)ze^{-(f^*(r)+g^*(r))z}
                                           \nonumber  \\ 
&  &  
+\frac{1}{2}e^{-\psi(f+g)}\int_E m(dr)\int n_I(dy)
(1-e^{-{f}(r)y})(1-e^{-{g}(r)y})               
                                             \\ 
& = & 
\frac{1}{2}e^{-\psi(f+g)} 
\left(\langle m,\frac{\alpha (f+g)^{\alpha-1}aK_B(f,g)}
{b+a(f+g)^{\alpha}}\rangle 
+\langle m,K_I(f,g)\rangle\right),           
\end{eqnarray*} 
which is finite. 
\end{lm}
{\it Proof.}~ 
It follows that 
\begin{eqnarray*} 
\Gamma(\Psi_f,\Psi_g) 
& = & 
\frac{1}{2}e^{-\lg \eta, f+g\rg}
\int_E \eta(dr)a(r)
\int n_B(dy)(1-e^{-f(r)y})(1-e^{-g(r)y})  \\ 
&  &  
+\frac{1}{2}e^{-\lg \eta, f+g\rg}
\int_E m(dr)\int n_I(dy)
(1-e^{-f(r)y})(1-e^{-{g}(r)y})            \\ 
& = & 
\frac{1}{2}e^{-\lg \eta, f+g\rg}
\left(\langle \eta,a K_B(f,g)\rangle 
+\langle m,K_I(f,g)\rangle\right).      
\end{eqnarray*} 
Note that the function $r\mapsto K_I(f(r),g(r))$ 
is an element of $B_+(E)$. 
Defining $h\in B_+(E)$ by 
$h(r)=a(r)K_B(f(r),g(r))$ and recalling that 
$\widetilde{\cE}(\Psi_f,\Psi_g)
=E^{Q_{\alpha}}[\Gamma(\Psi_f,\Psi_g)]$, 
we need only to show that 
\begin{eqnarray} 
I(f+g;h)
& := & 
E^{Q_{\alpha}}\left[e^{-\lg \eta, f+g\rg}
\lg \eta, h\rg \right]                       \nonumber \\ 
& = & 
e^{-\psi(f+g)}\int_E m_a(dr)
\frac{a(r)^{1/\alpha}}{b(r)^{1/\alpha}}h(r)
\int\Lambda(dz)ze^{-(f^*(r)+g^*(r))z} 
                                             \nonumber \\ 
& = & 
e^{-\psi(f+g)}
\langle m,\frac{\alpha (f+g)^{\alpha-1}h}
{b+a(f+g)^{\alpha}}\rangle 
                                                 \label{3.4}
\end{eqnarray} 
and that this is finite. 
The second equality can be verified to hold 
by (\ref{2.17}) and (\ref{3.1}) together: 
\begin{eqnarray*} 
I(f+g;h)  
& = & 
-\left.\frac{d}{d\epsilon}
E^{Q_{\alpha}}\left[e^{-\lg \eta,f+g+\epsilon h\rg}\right] 
\right|_{\epsilon=0} 
\ = \ 
-\left.\frac{d}{d\epsilon} 
e^{-\psi(f+g+\epsilon h)}\right|_{\epsilon=0}           \\ 
& = & 
e^{-\psi(f+g)}\int_E m_a(dr)h^*(r) 
\int\Lambda(dz)ze^{-(f^*(r)+g^*(r))z} \\ 
& = & 
e^{-\psi(f+g)}\int_E m_a(dr)
\frac{a(r)^{1/\alpha}}{b(r)^{1/\alpha}}h(r)
\int\Lambda(dz)ze^{-(f^*(r)+g^*(r))z}. 
\end{eqnarray*} 
For the proof of the last equality in (\ref{3.4}), 
we make use of another expression 
for $I(f+g;h)$ deduced from (\ref{2.17}) only: 
\begin{eqnarray*} 
I(f+g;h)  
& = & 
-\left.\frac{d}{d\epsilon} 
\exp\left[-\langle m,a^{-1}
\log(1+ab^{-1}(f+g+\epsilon h)^{\alpha}) 
\rangle\right]\right|_{\epsilon=0}           \\ 
& = & 
e^{-\psi(f+g)}
\langle m,\frac{\alpha (f+g)^{\alpha-1}h}
{b+a(f+g)^{\alpha}}\rangle. 
\end{eqnarray*} 
Here, by (\ref{3.2}) 
\[ 
0\le \alpha(f+g)^{\alpha-1}h
\le (f+g)^{\alpha-1}a(f+g)^{\alpha+1}=a(f+g)^{2\alpha} 
\] 
and so $I(f+g;h)$ is finite. 
\qed 

\medskip 

\noindent 
{\it Remark.}~ 
Noting that (\ref{3.4}) is clearly valid 
for every $h\in B_+(E)$ and 
combining (\ref{2.2}) with (\ref{3.4}), 
we get for any $f,g\in B_+(E)$ 
\begin{eqnarray*} 
E^{{Q}_{\alpha}}
\left[(-\cL_{\alpha})\Psi_f\cdot\Psi_g\right] 
& = & 
-E^{{Q}_{\alpha}}
\left[\Psi_{f+g}(\eta)\cdot
\frac{1}{\alpha}\langle \eta,af^{\alpha+1}+bf \rangle
-\Psi_{f+g}(\eta)\langle m,f^{\alpha} \rangle\right] \\ 
& = & 
-e^{-\psi(f+g)}\left(\langle m,
\frac{(f+g)^{\alpha-1}(af^{\alpha+1}+bf)}{b+a(f+g)^{\alpha}}
\rangle-\langle m,f^{\alpha} \rangle\right), 
\end{eqnarray*} 
from which the last expression in Lemma 3.1 
for the symmetric part 
\[ 
\widetilde{\cE}(\Psi_f,\Psi_g)
=
\frac{1}{2}\left(E^{{Q}_{\alpha}}
\left[(-\cL_{\alpha})\Psi_f\cdot\Psi_g\right]
+E^{{Q}_{\alpha}}
\left[(-\cL_{\alpha})\Psi_g\cdot\Psi_f\right]\right) 
\] 
can be recovered. 

\medskip 

Our objective is to show the positivity 
of ${\rm gap}(\overline{\cL_{\alpha}}^{(2)})$. 
The contribution here in this direction is the reduction 
to a certain estimate regarding the one-dimensional model. 
For a measurable function $f$ on $E$, 
the essential supremum (resp. the essential infimum) of $f$ 
with respect to $m$ is denoted by 
$\ds{\esssup f}$ (resp. $\ds{\essinf f}$). 
Let $D$ be the linear span of functions on $\R_+$ 
of the form $F_{\lambda}(z):=e^{-\lambda z}$ 
for some $\lambda> 0$.  
\begin{th}
Suppose that $b\in C_{++}(E)$. 
Let $\gamma>0$ be a constant. 
If for every $F\in D$ 
\begin{eqnarray} 
\lefteqn{\int\Lambda(dz)(F(z)-F(0))^2}        \nonumber \\ 
& \le & 
\frac{\gamma}{2}
\left[
\int\Lambda(dz)z\int n_B(dy)(F(z+y)-F(z))^2
+\int n_I(dy)(F(y)-F(0))^2\right], 
                                               \label{3.5}  
\end{eqnarray} 
then for any $\Psi\in \cD$  
\be 
{\rm var}(\Psi) \le 
\gamma \esssup(b^{-1})\cE(\Psi) 
                                               \label{3.6}  
\ee 
and it holds that 
$\ds{{\rm gap}(\overline{\cL_{\alpha}}^{(2)})
\ge \gamma^{-1}\essinf b}$. 
\end{th}
This kind of reduction was discovered 
by Stannat \cite{Stannat05} (Theorem 1.2) 
for a lower estimate for the quadratic form of gradient type. 
In particular, for the process associated with $\cL_1$ 
in Remark at the end of Section 2, 
the condition corresponding to (\ref{3.5}) reads 
\be 
\int \Lambda_1(dz)(F(z)-F(0))^2 
\le \gamma\int \Lambda_1(dz)z(F'(z))^2,   
                                               \label{3.7} 
\ee 
where $\Lambda_1(dz)=z^{-1}e^{-z}dz$ is 
the L\'evy measure of a gamma distribution. 
While (\ref{3.7}) with $\gamma=1$ is verified easily 
by applying Schwarz's inequality to 
$F(z)-F(0)=\int_0^zF'(w)dw$,  
showing an inequality of the form (\ref{3.5}) is 
more difficult and we postpone it until the next section. 
However, as will be seen below, 
the reduction itself is proved 
in a similar way to \cite{Stannat05}. \\ 
{\it Proof of Theorem 3.2.}~ 
Consider a function $\Psi$ expressed as 
a finite sum $\Psi=\sum_i c_i \Psi_{f_i}$, where   
$c_i\in\R$ and $f_i\in B_{+}(E)$. 
Putting $d_i=c_ie^{-\psi(f_i)}$, 
observe from (\ref{3.1}) that 
\begin{eqnarray} 
{\rm var}(\Psi) 
& = & 
\sum_{i,j}c_ic_j
\left(e^{-\psi(f_i+f_j)}
-e^{-\psi(f_i)}e^{-\psi(f_j)}\right)  
                                             \nonumber \\ 
& = & 
\sum_{i,j}d_id_j
\left(e^{\psi(f_i)+\psi(f_j)-\psi(f_i+f_j)}-1\right)  
                                             \nonumber \\ 
& = & 
\sum_{i,j}d_id_j
\left[\exp\left(\int_E\int m_a(dr)\Lambda(dz)
(1-e^{-f^*_i(r)z})
(1-e^{-f^*_j(r)z})\right)-1\right]            
                                             \nonumber \\ 
& = & 
\sum_{i,j}d_id_j
\sum_{N=1}^{\infty}\frac{1}{N!}
\left(\int_E\int m_a(dr)\Lambda(dz)
(1-e^{-f^*_i(r)z})
(1-e^{-f^*_j(r)z})\right)^N.           \label{3.8}
\end{eqnarray} 
Rewrite in terms of the 
$N$-fold product measures 
$m_a^{\otimes N}$ and $\Lambda^{\otimes N}$ to 
obtain the following disintegration formula 
for the variance functional:  
\begin{eqnarray} 
{\rm var}(\Psi) 
& = & 
\sum_{N=1}^{\infty}\frac{1}{N!}\sum_{i,j}d_id_j
\int_{E^N}\int_{\R_+^N} 
m_a^{\otimes N}(d\br_N)\Lambda^{\otimes N}(d\bz_N)   
                                             \nonumber \\ 
&  & 
\ \prod_{k=1}^N(1-e^{-f^*_i(r_k)z_k})
\prod_{l=1}^N(1-e^{-f^*_j(r_l)z_l})    \nonumber \\ 
& = & 
\sum_{N=1}^{\infty}\frac{1}{N!}
\int_{E^N}\int_{\R_+^N} 
m_a^{\otimes N}(d\br_N)\Lambda^{\otimes N}(d\bz_N)
\left[\sum_{i}d_i\prod_{k=1}^N
(1-e^{-f^*_i(r_k)z_k})\right]^2,    \label{3.9}
\end{eqnarray} 
where $\br_N=(r_1,\ldots,r_N)$ and 
$\bz_N=(z_1,\ldots,z_N)$. 
Given $\br_N=(r_1,\ldots,r_N)\in E^N$ and 
$z_1,\ldots,z_{N-1}\in\R_+$ arbitrarily, 
apply (\ref{3.5}) to the function 
\[ 
z_N\mapsto \sum_{i}d_i\left\{\prod_{k=1}^{N-1}
(1-e^{-f^*_i(r_k)z_k})\right\}
e^{-f^*_i(r_N)z_N}
\] 
to get 
\begin{eqnarray} 
\lefteqn{
\frac{2}{\gamma}\int \Lambda(dz_N)
\left[\sum_{i}d_i\prod_{k=1}^N
(1-e^{-f^*_i(r_k)z_k})\right]^2}  \nonumber \\ 
& \le & 
\int\Lambda(dz)z\int n_B(dy)
\left[\sum_{i}d_i\prod_{k=1}^{N-1}
(1-e^{-f^*_i(r_k)z_k})
(e^{-f^*_i(r_N)(z+y)}
-e^{-f^*_i(r_N)z})\right]^2
                                       \nonumber \\ 
&  & 
+\int n_I(dy)\left[\sum_{i}d_i\prod_{k=1}^{N-1}
(1-e^{-f^*_i(r_k)z_k})
(e^{-f^*_i(r_N)y}-1)\right]^2    \nonumber \\
& = & 
\frac{a(r_N)^{1+1/\alpha}}{b(r_N)^{1+1/\alpha}}
\int\Lambda(dz)z
\int n_B(dy)
\left[\sum_{i}d_i\prod_{k=1}^{N-1}
(1-e^{-f^*_i(r_k)z_k})e^{-f^*_i(r_N)z}
(1-e^{-{f_i}(r_N)y})\right]^2
                                       \nonumber \\ 
&  & 
+\frac{a(r_N)}{b(r_N)}
\int n_I(dy)\left[\sum_{i}d_i\prod_{k=1}^{N-1}
(1-e^{-f^*_i(r_k)z_k})
(1-e^{-f_i(r_N)y})\right]^2. 
                                           \label{3.10}
\end{eqnarray} 
Here, a suitable change of variable has been made 
for each integral with respect to $n_B(dy)$ and $n_I(dy)$ 
in order to replace $f^*_i(r_N)y$ by ${f_i}(r_N)y$.

Set $C=\ds{\esssup(b^{-1})}$ so that 
\[ 
m_a(dr_N)
=a(r_N)^{-1}m(dr_N)
\le C\cdot 
b(r_N)a(r_N)^{-1}m(dr_N)
\] 
in distributional sense. 
Combining (\ref{3.9}) with (\ref{3.10}), 
we can dominate $2{\rm var}(\Psi)/\gamma$ by 
\begin{eqnarray*} 
\lefteqn{
\sum_{N=1}^{\infty}\frac{C}{N!}
\int_{E^{N-1}}\int_{\R_+^{N-1}} 
m_a^{\otimes {N-1}}(d\br_{N-1})
\Lambda^{\otimes {N-1}}(d\bz_{N-1})
\int_E m(dr)\frac{a(r)^{1/\alpha}}{b(r)^{1/\alpha}}}  \\ 
& & 
\int\Lambda(dz)z\int n_B(dy)
\left[\sum_{i}d_i\prod_{k=1}^{N-1}
(1-e^{-f^*_i(r_k)z_k})e^{-f^*_i(r)z}
(1-e^{-{f_i}(r)y})\right]^2
                                       \\ 
&  &  
+\sum_{N=1}^{\infty}\frac{C}{N!}
\int_{E^{N-1}}\int_{\R_+^{N-1}} 
m_a^{\otimes {N-1}}(d\br_{N-1})
\Lambda^{\otimes {N-1}}(d\bz_{N-1})\int_E m(dr) \\ 
&  & 
\int n_I(dy)\left[\sum_{i}d_i\prod_{k=1}^{N-1}
(1-e^{-f^*_i(r_k)z_k})(1-e^{-f_i(r)y})\right]^2   \\ 
& \le &  
\sum_{N=0}^{\infty}\frac{C}{N!}
\int_{E^{N}}\int_{\R_+^{N}}m_a^{\otimes {N}}(d\br_{N})
\Lambda^{\otimes {N}}(d\bz_{N})
\int_E m(dr)\frac{a(r)^{1/\alpha}}{b(r)^{1/\alpha}}  \\ 
&   & 
\int\Lambda(dz)z\int n_B(dy)
\left[\sum_{i}d_i\prod_{k=1}^{N}
(1-e^{-f^*_i(r_k)z_k})e^{-f^*_i(r)z}
(1-e^{-{f_i}(r)y})\right]^2
                                       \\ 
&  &  
+\sum_{N=0}^{\infty}\frac{C}{N!}
\int_{E^{N}}\int_{\R_+^{N}} 
m_a^{\otimes {N}}(d\br_{N})
\Lambda^{\otimes {N}}(d\bz_{N})\int_E m(dr) \\ 
&  & 
\int n_I(dy)\left[\sum_{i}d_i\prod_{k=1}^{N}
(1-e^{-f^*_i(r_k)z_k})
(1-e^{-f_i(r)y})\right]^2           \\ 
& = &  
\sum_{i,j}d_id_j
\sum_{N=0}^{\infty}\frac{C}{N!}
\left(\int_{E}\int m_a(dr_1)\Lambda(dz_1)
(1-e^{-f^*_i(r_1)z_1})
(1-e^{-f^*_j(r_1)z_1})\right)^N               \\ 
&   & 
\int_E m(dr)\frac{a(r)^{1/\alpha}}{b(r)^{1/\alpha}}
\int n_B(dy)(1-e^{-{f_i}(r)y})(1-e^{-{f_j}(r)y})
\int\Lambda(dz)ze^{-(f^*_i(r)+f^*_j(r))z}
                                       \\ 
&  &  
+\sum_{i,j}d_id_j
\sum_{N=0}^{\infty}\frac{C}{N!}
\left(\int_{E}\int m_a(dr_1)\Lambda(dz_1)
(1-e^{-f^*_i(r_1)z_1})
(1-e^{-f^*_j(r_1)z_1})\right)^N  \\ 
& & 
\int_E m(dr)\int n_I(dy)
(1-e^{-{f_i}(r)y})(1-e^{-{f_j}(r)y})           \\ 
& = &  
C\sum_{i,j}c_ic_j
e^{-\psi(f_i+f_j)}\int_E m(dr)
\frac{a(r)^{1/\alpha}}{b(r)^{1/\alpha}}
\int n_B(dy)(1-e^{-{f_i}(r)y})(1-e^{-{f_j}(r)y})  \\ 
&   & 
\int\Lambda(dz)ze^{-(f^*_i(r)+f^*_j(r))z}
                                              \\ 
&  &  
+C\sum_{i,j}c_ic_j
e^{-\psi(f_i+f_j)}\int_E m(dr)\int n_I(dy)
(1-e^{-{f_i}(r)y})(1-e^{-{f_j}(r)y}),   
\end{eqnarray*} 
where the last two equalities are seen by 
similar calculations to (\ref{3.8}) and (\ref{3.9}). 
Since the symmetric part $\widetilde{\cE}$ 
of Dirichlet form is bilinear, 
(\ref{3.6}) for $\Psi\in\cD$ follows from Lemma 3.1.

It remains to prove that (\ref{3.6}) extends to 
$\Psi\in D(\overline{\cL_{\alpha}}^{(2)})$. 
Since $(\overline{\cL_{\alpha}}^{(2)},
D(\overline{\cL_{\alpha}}^{(2)}))$ is the closure of 
$(\cL_{\alpha},\cF)$ in $L^2(Q_{\alpha})$, 
we need only to show that 
(\ref{3.6}) extends to $\Psi\in \cF$. 
Given $\Psi\in \cF$, we see from Theorem 2.2 (i) 
that there exists a sequence 
$\{\Psi_{N}\}_{N=1}^{\infty}\subset \cD$ such that 
\[ 
\Vert \Psi_N-\Psi\Vert_{\infty}+
\Vert \cL_{\alpha}\Psi_N-\cL_{\alpha}\Psi\Vert_{\infty}
\to 0 \quad \mbox{as} \quad N\to\infty. 
\] 
Hence  
\[ 
\Vert \Psi_N-\Psi\Vert_{L^2(Q_{\alpha})}+\cE(\Psi_N-\Psi)
\to 0 \quad \mbox{as} \quad N\to\infty.  
\] 
This implies that (\ref{3.6}) holds for any $\Psi\in \cF$ 
and we complete the proof of Theorem 3.2. \qed 

\medskip 

\noindent
{\it Remark.}~ In view of calculations in the proof of 
Theorem 3.2, it is clear that under the same assumption 
an appropriate version of 
the inequality (\ref{3.6}) 
holds for a more general class of functions $a$ and $b$. 
To be more precise, suppose that $a,b\in B_+(E)$ are 
uniformly positive and that $Q_{\alpha}$ has 
Laplace transform (\ref{2.17}). Then, 
assuming that (\ref{3.5}) is valid for all $F\in D$, 
we have 
\[ 
E^{Q_{\alpha}}
\left[\left(\Psi-E^{Q_{\alpha}}[\Psi]\right)^2\right] 
\le 
\gamma \esssup(b^{-1})
E^{Q_{\alpha}}\left[\Gamma(\Psi,\Psi)\right]  
\] 
for any $\Psi\in\cD$, where 
\begin{eqnarray*} 
\Gamma(\Psi,\Psi)(\eta) 
& = & 
\frac{1}{2}\int n_B(dz)\int_E\eta(dr)a(r)
\left[\Psi(\eta+z\delta_r)-\Psi(\eta)\right]^2 \\ 
& & 
+\frac{1}{2}\int n_I(dz)\int_Em(dr)
\left[\Psi(\eta+z\delta_r)-\Psi(\eta)\right]^2. 
\end{eqnarray*} 
This fact reflects the convolution property with respect to $m$ 
mentioned in the Introduction. 
(Note that the condition (\ref{3.5}) is independent of $m$.) 

\section{Spectral gap for the $\alpha$-CIR model}
\setcounter{equation}{0} 
This section is devoted to the proof of (\ref{3.5}) 
for some $0<\gamma<\infty$. The strategy should be 
different from the one already mentioned for (\ref{3.7}) 
with $\Lambda_1(dz)=z^{-1}e^{-z}dz$ 
at least because no informative expression for 
the density of $\Lambda$ in (\ref{3.1}) appears to be available. 
Let us illustrate another approach 
we will take and call `the method of intrinsic kernel' 
by revisiting (\ref{3.7}). 
Suppose that $F\in D$ is a finite sum 
$F=\sum_ic_iF_{\lambda_i}$. 
We will use the notation $\one_S$ 
standing for the indicator function of a set $S$ 
and $\partial_t=\partial/\partial t$ for simplicity. 
Letting 
$\psi_1(\lambda)=\log(1+\lambda)
=\int\Lambda_1(dz)(1-e^{-\lambda z})$, 
observe that 
\begin{eqnarray} 
{\bf U}_1(F) \ := \ 
\int\Lambda_1(dz)(F(z)-F(0))^2  
& = & 
\sum_{i,j}c_ic_j
\int\Lambda_1(dz)(1-e^{-\lambda_i z})(1-e^{-\lambda_j z}) 
                                         \nonumber \\ 
& = & 
\sum_{i,j}c_ic_j
\left(-\psi_1(\lambda_i+\lambda_j)
+\psi_1(\lambda_i)+\psi_1(\lambda_j)\right)  \nonumber \\
& = & 
\sum_{i,j}c_ic_j\int_0^{\lambda_i}ds
\int_0^{\lambda_j}dt(-\psi_1''(s+t))     \nonumber \\
& = & 
\int ds\int dt \overline{F}(s)\overline{F}(t) 
(-\psi_1''(s+t)),   
                                                \label{4.1}
\end{eqnarray} 
where $\overline{F}(s)=\sum_i c_i\one_{[0,\lambda_i]}(s)$. 
On the other hand, by putting $K_1(s,t)=st\psi_1'(s+t)$ 
\begin{eqnarray} 
{\bf V}_1(F) \ := \ 
\int\Lambda_1(dz)z(F'(z))^2  
& = & 
\sum_{i,j}c_ic_j\lambda_i\lambda_j
\int\Lambda_1(dz)ze^{-\lambda_i z}e^{-\lambda_j z} 
                                         \nonumber \\ 
& = & 
\sum_{i,j}c_ic_jK_1(\lambda_i,\lambda_j) \nonumber \\
& = & 
\sum_{i,j}c_ic_j\int_0^{\lambda_i}ds
\int_0^{\lambda_j}dt 
\partial_s\partial_t K_1(s,t) \nonumber \\
& = & 
\int ds\int dt \overline{F}(s)\overline{F}(t) 
\partial_s\partial_t K_1(s,t).  
                                               \label{4.2}
\end{eqnarray} 
It is reasonable to call  
$\partial_s\partial_t K_1$ 
the intrinsic kernel of the quadratic form ${\bf V}_1$. 
Similarly, (\ref{4.1}) shows that 
the intrinsic kernel of ${\bf U}_1$ is 
the function $(s,t)\mapsto -\psi_1''(s+t)$.

Given two symmetric measurable functions 
$J$ and $K$ on $\R_+^2$, we write $K\gg J$ 
if $K-J$ is nonnegative definite in the sense that 
\[ 
\int ds\int dt G(s)G(t) 
(K(s,t)-J(s,t)) \ge 0 
\] 
for any bounded Borel function $G$ 
on $\R_+$ with compact support. 
By virtue of Fubini's theorem, 
$K\gg 0$ if $K$ is of `canonical form' 
\[ 
K(s,t)=\int_SM(d\omega)\sigma(s,\omega)\sigma(t,\omega), 
\quad s,t\in\R_+ 
\] 
for some measure space $(S,M)$ and 
measurable function $\sigma$ on $\R_+\times S$. 
In view of (\ref{4.1}) and (\ref{4.2}), 
it is clear that the inequality 
$\gamma {\bf V}_1(F)\ge {\bf U}_1(F)$ is implied by 
\[ 
\gamma \partial_s\partial_t K_1(s,t)
+\psi_1''(s+t) \gg 0. 
\]  
For $\gamma=1$, 
this holds true since by direct calculations 
\[ 
\partial_s\partial_t K_1(s,t)+\psi_1''(s+t)                                    
\ = \  
\frac{2st}{(1+s+t)^3}                      \\ 
\ = \
\int dz z^2e^{-z}se^{-sz}te^{-tz}, 
\] 
which is of canonical form. Furthermore, 
this expression makes it possible to identify 
the associated `remainder form': 
\begin{eqnarray*} 
{\bf V}_1(F)-{\bf U}_1(F)  
& = & 
\int ds\int dt \overline{F}(s)\overline{F}(t) 
\int dz z^2e^{-z}se^{-sz}te^{-tz}             \\ 
& = & 
\int dz e^{-z}
\left(\int ds \overline{F}(s)sze^{-sz}\right)^2   \\ 
& = & 
\int dz e^{-z}
\left(\sum_ic_i\int_0^{\lambda_i} ds sze^{-sz}\right)^2 
                                                 \\ 
& = & 
\int dz e^{-z}\left(\sum_ic_i
\frac{e^{-\lambda_iz}-1+\lambda_ize^{-\lambda_iz}}{z}\right)^2
                                                 \\ 
& = & 
\int\Lambda_1(dz)z^{-1}(F(z)-F(0)-zF'(z))^2.   
\end{eqnarray*}
It should be emphasized that the above calculations 
require only an explicit form of 
the Laplace exponent $\psi_1$. 

Turning to the case $0<\alpha<1$, 
we adopt the method of intrinsic kernels 
to show (\ref{3.5}) for $\Lambda$ such that 
\be 
\psi(\lambda):=\log(1+\lambda^{\alpha})
=\int\Lambda(dz)(1-e^{-\lambda z}), 
\quad \lambda\ge 0. 
                                   \label{4.3} 
\ee 
(We continue to adopt this notation 
as it is a one-dimensional version of (\ref{3.1}).) 
Namely, we shall  \\ 
(I) calculate the intrinsic kernels of  
${\bf U}(F):=\int\Lambda(dz)(F(z)-F(0))^2$ 
and of 
\[ 
{\bf V}(F):=
\frac{1}{2}\left[
\int\Lambda(dz)z\int n_B(dy)(F(z+y)-F(z))^2
+\int n_I(dy)(F(y)-F(0))^2\right], 
\] 
and then \\ 
(II) compare the two kernels 
as nonnegative definite functions. \\ 
The following lemma concerns the step (I). 
\begin{lm}
The intrinsic kernels of ${\bf U}$ and of ${\bf V}$ 
are respectively given by $J(s,t):=-\psi''(s+t)$ and 
$\widetilde{K}(s,t)
:=\partial_s\partial_t K(s,t)$, where 
\begin{eqnarray} 
2K(s,t)
& = & 
\frac{(s+t)^{\alpha-1}}{1+(s+t)^{\alpha}}
\left[(s+t)^{\alpha+1}-s^{\alpha+1}-t^{\alpha+1}\right]
+\left[s^{\alpha}+t^{\alpha}-(s+t)^{\alpha}\right].  
                                            \label{4.4} 
\end{eqnarray} 
\end{lm}
{\it Proof.}~
The intrinsic kernel of ${\bf U}$ is deduced 
in the same manner as (\ref{4.1}). 
The derivation of (\ref{4.4}) is similar to 
the calculations at the beginning of the proof of Lemma 3.1. 
Indeed, for $F=\sum_ic_iF_{\lambda_i}$ 
\begin{eqnarray*} 
{\bf V}(F) 
& = & 
\frac{1}{2}\sum_{i,j}c_ic_j
\int\Lambda(dz)ze^{-(\lambda_i+\lambda_j) z} 
\int n_B(dy)(1-e^{-\lambda_i y})(1-e^{-\lambda_j y})  \\ 
&  & 
+\frac{1}{2}\sum_{i,j}c_ic_j\int n_I(dy)
(1-e^{-\lambda_i y})(1-e^{-\lambda_j y}) \\ 
& = & 
\sum_{i,j}c_ic_jK(\lambda_i,\lambda_j), 
\end{eqnarray*} 
where the last equality is seen from 
$\int\Lambda(dz)ze^{-\lambda z}=\psi'(\lambda)$, 
(\ref{4.3}), (\ref{3.2}) and (\ref{3.3}).  
The rest of the proof is 
the same as (\ref{4.2}) with $K$ in place of $K_1$. \qed 

\medskip 

\noindent 
Remark that, as $\alpha \uparrow 1$, the right side 
of (\ref{4.4}) tends to $\frac{2st}{1+s+t}=2K_1(s,t)$. 
The main result of this section is 
obtained by accomplishing not only the step (II) 
but also identification of the remainder form. 
\begin{th}
For each $F\in D$, 
\begin{eqnarray} 
2{\bf V}(F)-{\bf U}(F)
& = & 
\frac{\alpha}{\alpha+1}
\int\Lambda(dz)z^3\int n_B(dy)  
\left[\frac{F(z)-F(0)}{z}
-\frac{F(z+y)-F(0)}{z+y}\right]^2   \nonumber \\ 
&  & 
+\int\Lambda(dz)z^2\int n_I(dy) 
\left[\frac{F(z)-F(0)}{z}
-\frac{F(z+y)-F(0)}{z+y}\right]^2   \nonumber \\ 
&  & 
+ \frac{1}{\alpha+1}\int\Lambda(dz)z\int n_B(dy)
\left(F(z+y)-F(z)\right)^2.                 \label{4.5} 
\end{eqnarray} 
In particular, (\ref{3.5}) with $\gamma=2$ holds true. 
\end{th}
{\it Proof.}~
Let $J,K$ and $\widetilde{K}$ be as in Lemma 4.1. 
Recalling (\ref{3.2}) and (\ref{3.3}), 
we will exploit the following expression 
for $2K$ in (\ref{4.4}): 
\be  
2K(s,t)=\rho(s+t)K_B(s,t)+K_I(s,t), 
                                            \label{4.6} 
\ee 
where 
\[  
\rho(\lambda)
=\frac{\alpha\lambda^{\alpha-1}}{1+\lambda^{\alpha}}
=\psi'(\lambda)
=\int\Lambda(dz)ze^{-\lambda z}. 
\] 
Differentiating (\ref{4.6}) in $s$ and $t$ yields 
\begin{eqnarray}
2\widetilde{K}(s,t)
& = & 
\partial_s\partial_t\left[\rho(s+t)K_B(s,t)\right]
+\alpha(1-\alpha)(s+t)^{\alpha-2}         \label{4.7}  \\ 
& = & 
\rho''(s+t)K_B(s,t)
+\rho'(s+t)(\partial_s+\partial_t)K_B(s,t) \nonumber \\ 
& & 
+\rho(s+t)(\alpha+1)(s+t)^{\alpha-1}
+\alpha(1-\alpha)(s+t)^{\alpha-2}. 
                                            \label{4.8} 
\end{eqnarray}
Since 
\[ 
(\partial_s+\partial_t)K_B(s,t) 
=\frac{\alpha+1}{\alpha}
\left[-K_I(s,t)+(s+t)^{\alpha}\right], 
\] 
(\ref{4.8}) becomes 
\begin{eqnarray}
2\widetilde{K}(s,t)
& = & 
\rho''(s+t)K_B(s,t)
-\frac{\alpha+1}{\alpha}\rho'(s+t)K_I(s,t)
+\frac{\alpha+1}{\alpha}\rho'(s+t)(s+t)^{\alpha} \nonumber \\ 
& & 
+\rho(s+t)(\alpha+1)(s+t)^{\alpha-1}
+\alpha(1-\alpha)(s+t)^{\alpha-2}. 
                                            \label{4.9} 
\end{eqnarray}
Here, it is direct to see that 
the sum of the last three terms on the right side equals 
\[ 
-\frac{\alpha+1}{\alpha}\rho'(s+t)-(1-\alpha)(s+t)^{\alpha-2} 
\] 
or equivalently that 
\[ 
\frac{1}{\alpha}\rho'(\lambda)(\lambda^{\alpha}+1)
+\rho(\lambda)\lambda^{\alpha-1}
=-(1-\alpha)\lambda^{\alpha-2}. 
\] 
By the above identity and (\ref{4.9}) together 
\begin{eqnarray*}
2\widetilde{K}(s,t)
& = & 
\rho''(s+t)K_B(s,t)
-\frac{\alpha+1}{\alpha}\rho'(s+t)K_I(s,t)         \\ 
& & 
-\frac{\alpha+1}{\alpha}
\rho'(s+t)-(1-\alpha)(s+t)^{\alpha-2} \\ 
& = & 
\rho''(s+t)K_B(s,t)
-\frac{\alpha+1}{\alpha}\rho'(s+t)K_I(s,t)          \\ 
& & 
+\frac{\alpha+1}{\alpha}J(s,t)
-(1-\alpha)(s+t)^{\alpha-2}             \\ 
& = & 
\rho''(s+t)K_B(s,t)
-\frac{\alpha+1}{\alpha}\rho'(s+t)K_I(s,t)        \\ 
& & 
+\frac{\alpha+1}{\alpha}J(s,t)
-\frac{2}{\alpha}\widetilde{K}(s,t)
+\frac{1}{\alpha}
\partial_s\partial_t\left[\rho(s+t)K_B(s,t)\right], 
\end{eqnarray*} 
where the last equality follows from (\ref{4.7}). 
Consequently 
\begin{eqnarray*}
\lefteqn{2\widetilde{K}(s,t)-J(s,t)}          \\ 
& = & 
\frac{\alpha}{\alpha+1}\rho''(s+t)K_B(s,t)
-\rho'(s+t)K_I(s,t) 
+\frac{1}{\alpha+1}
\partial_s\partial_t\left[\rho(s+t)K_B(s,t)\right].      
\end{eqnarray*} 
Each of the terms on the right side is 
nonnegative definite because of 
\[ 
\rho''(s+t)K_B(s,t)
=\int\Lambda(dz)z^3
\int n_B(dy)e^{-sz}(1-e^{-sy})e^{-tz}(1-e^{-ty}) ,   
\] 
\[ 
-\rho'(s+t)K_I(s,t)
=\int\Lambda(dz)z^2
\int n_I(dy)e^{-sz}(1-e^{-sy})e^{-tz}(1-e^{-ty})    
\] 
and 
\[ 
\partial_s\partial_t
\left[\rho(s+t)K_B(s,t)\right]
=\int\Lambda(dz)z\int n_B(dy)
\partial_s\left[e^{-sz}(1-e^{-sy})\right]
\partial_t\left[e^{-tz}(1-e^{-ty})\right].      
\] 
Therefore, $2\widetilde{K}\gg J$ and 
so $2{\bf V}(F)\ge {\bf U}(F)$ for any $F\in D$. 
We further proceed to identify the remainder form 
$F\mapsto 2{\bf V}(F)-{\bf U}(F)$. 
With the help of the canonical representations 
in the above, we deduce 
$2{\bf V}(F)-{\bf U}(F)
=\sum_{i,j}c_ic_jK'(\lambda_i,\lambda_j)$ 
for $F=\sum_ic_iF_{\lambda_i}$, 
where 
\begin{eqnarray*} 
\lefteqn{K'(\lambda_i,\lambda_j) \ = \ 
\int_0^{\lambda_i}ds\int_0^{\lambda_j}dt
\left(2\widetilde{K}(s,t)-J(s,t)\right)}                  \\ 
& = & 
\frac{\alpha}{\alpha+1}\int\Lambda(dz)z^3\int n_B(dy)
\left[\frac{1-e^{-\lambda_iz}}{z}
-\frac{1-e^{-\lambda_i(z+y)}}{z+y}\right] 
\left[\frac{1-e^{-\lambda_jz}}{z}
-\frac{1-e^{-\lambda_j(z+y)}}{z+y}\right]       \\ 
&  & 
+\int\Lambda(dz)z^2\int n_I(dy)
\left[\frac{1-e^{-\lambda_iz}}{z}
-\frac{1-e^{-\lambda_i(z+y)}}{z+y}\right] 
\left[\frac{1-e^{-\lambda_jz}}{z}
-\frac{1-e^{-\lambda_j(z+y)}}{z+y}\right]                  \\ 
&  & 
+\frac{1}{\alpha+1}\int\Lambda(dz)z\int n_B(dy)
\left[e^{-\lambda_iz}(1-e^{-\lambda_iy})\right]
\left[e^{-\lambda_jz}(1-e^{-\lambda_jy})\right]. 
\end{eqnarray*} 
Accordingly 
\begin{eqnarray*} 
2{\bf V}(F)-{\bf U}(F)
& = & 
\frac{\alpha}{\alpha+1}
\int\Lambda(dz)z^3\int n_B(dy)  
\left[\sum_{i}c_i\left(\frac{1-e^{-\lambda_iz}}{z}
-\frac{1-e^{-\lambda_i(z+y)}}{z+y}\right)\right]^2   \\ 
&  & 
+\int\Lambda(dz)z^2\int n_I(dy) 
\left[\sum_{i}c_i\left(\frac{1-e^{-\lambda_iz}}{z}
-\frac{1-e^{-\lambda_i(z+y)}}{z+y}\right)\right]^2   \\ 
&  & 
+ \frac{1}{\alpha+1}\int\Lambda(dz)z\int n_B(dy)
\left[\sum_{i}c_i\left(e^{-\lambda_iz}
-e^{-\lambda_i(z+y)}\right)\right]^2. 
\end{eqnarray*} 
This coincides with the right side of 
(\ref{4.5}) for $F=\sum_ic_iF_{\lambda_i}$, 
and the proof of Theorem 4.2 is complete. 
\qed 

\medskip 

\noindent 
Combining Theorem 4.2 with Theorem 3.2 
gives immediately a lower estimate 
\[ 
{\rm gap}(\overline{\cL_{\alpha}}^{(2)}) 
\ge \frac{1}{2}\essinf b 
\] 
for ergodic measure-valued $\alpha$-CIR models 
discussed in the previous section. 
In fact, this bound is optimal as seen in 
\begin{th}
Let $\cL_{\alpha}$ be of the form (\ref{2.1}) 
for some $a,b\in C_{++}(E)$ and $m\in\cM(E)^{\circ}$. 
Then 
\[  
{\rm gap}(\overline{\cL_{\alpha}}^{(2)}) 
= \frac{1}{2}\essinf b. 
\] 
\end{th}
{\it Proof.}~
For the aforementioned reason, 
it suffices to show the upper estimate 
\[ 
{\rm gap}(\overline{\cL_{\alpha}}^{(2)}) 
\le \frac{1}{2}\essinf b. 
\] 
To this end, we use (a variant of) 
a characterization due to Liggett (\cite{Liggett89}, (2.5)): 
\[ 
{\rm gap}(\overline{\cL_{\alpha}}^{(2)})
=\inf_{t>0}\frac{1}{2t}
\inf\left\{-\log \frac{{\rm var}(T^2(t)\Psi)}{{\rm var}(\Psi)}
: 0<{\rm var}(\Psi)<\infty \right\}. 
\] 
This implies that 
\be 
{\rm gap}(\overline{\cL_{\alpha}}^{(2)})
\le \liminf_{t\to\infty}
\frac{1}{2t}(-\log{\rm var}(T^2(t)\Psi)) 
                                                 \label{4.10} 
\ee 
for any $\Psi$ such that $0<{\rm var}(\Psi)<\infty$. 
We now take $\Psi(\eta)=\exp(-\eta(E))$, 
for which $T^2(t)\Psi(\eta)$ is given by 
the right side of (\ref{2.15}) with $f=\one_E=:\one$.  
Recalling that the log-Laplace functional of $Q_{\alpha}$ 
is $\psi(f)=\lg m_a, \log(1+ab^{-1}f^{\alpha})\rg$, 
one can derive by (\ref{2.15}) 
\begin{eqnarray*} 
{\rm var}(T^2(t)\Psi)
& = & 
E^{Q_{\alpha}}\left[\exp\left(-2\lg \eta, V_t\one\rg
-2\int_0^t\lg m, (V_s\one)^{\alpha}\rg ds\right)\right]
-e^{-2\psi(\one)}                                        \\ 
& = & 
e^{-2\psi(\one)}\left(\exp\left[-\psi(2V_t\one)
+2\psi(V_t\one)\rg\right]-1\right),  
\end{eqnarray*} 
where the last equality is deduced from 
\[ 
\psi(\one)-\int_0^t\lg m, (V_s\one)^{\alpha}\rg ds
=\psi(V_t\one). 
\] 
Since by (\ref{2.16}) 
$(V_t\one(r))^{\alpha}\le e^{-tb(r)}$,  
$\Delta(t):=2\psi(V_t\one)-\psi(2V_t\one)\to 0$ 
as $t\to \infty$ and so 
\[ 
\liminf_{t\to\infty}
\frac{1}{t}(-\log{\rm var}(T^2(t)\Psi)) 
=
-\limsup_{t\to\infty}
\frac{1}{t}\log\Delta(t).    
\] 
By virtue of (\ref{4.10}), the proof of 
Theorem 4.3 reduces to showing that 
\be 
-\limsup_{t\to\infty}\frac{1}{t}\log\Delta(t)
\le \essinf b. 
                                               \label{4.11}  
\ee 
By straightforward calculations 
\begin{eqnarray} 
\Delta(t) 
& = &  
\lg m_a, \log
\left(1+\frac{\frac{a}{b}(2-2^{\alpha})(V_t\one)^{\alpha}
+\left(\frac{a}{b}\right)^2(V_t\one)^{2\alpha}}
{1+\frac{a}{b}2^{\alpha}(V_t\one)^{\alpha}}\right)\rg   
                                            \label{4.12} \\ 
& \ge & 
\lg m_a, \log
\left(1+\frac{\frac{a}{b}(2-2^{\alpha})(V_t\one)^{\alpha}}
{1+\frac{a}{b}2^{\alpha}(V_t\one)^{\alpha}}\right)\rg. \nonumber 
\end{eqnarray} 
Further, with the help of the inequality 
$\log(1+z)\ge z/(1+z)$ for $z\ge0$, we get 
\[ 
\Delta(t) 
\ge 
\lg m_a, 
\frac{\frac{a}{b}(2-2^{\alpha})(V_t\one)^{\alpha}}
{1+\frac{a}{b}2(V_t\one)^{\alpha}}\rg 
\ge 
\lg m, (V_t\one)^{\alpha}\rg 
\essinf \frac{2-2^{\alpha}}{b+2a}. 
\] 
Here, again by (\ref{2.16}) 
\[ 
(V_t\one(r))^{\alpha}
\ge e^{-tb(r)}\essinf \frac{b}{b+a}, 
\] 
and therefore 
\[ 
\limsup_{t\to\infty}
\frac{1}{t}\log\Delta(t) 
\ge   
\lim_{t\to\infty}
\frac{1}{t}\log\lg m, e^{-tb}\rg  
= 
\esssup(-b)\ =\ -\essinf b. 
\] 
This establishes (\ref{4.11}) and 
completes the proof of Theorem 4.3. \qed 

\medskip 

\noindent 
{\it Remarks.}~
(i) The same argument may apply to the case $\alpha=1$. 
But the resulting bound exhibits discontinuity at $\alpha=1$. 
Indeed, for $\alpha=1$ (\ref{4.12}) becomes 
\[ 
\Delta(t) 
=
\lg m_a, \log
\left(1+\frac{\left(\frac{a}{b}\right)^2(V_t\one)^{2}}
{1+\frac{a}{b}2V_t\one}\right)\rg,  
\] 
where $V_t\one(r)$ is given by the right side of 
(\ref{2.16}) with $\alpha=1$ and $f(r)=1$. 
(See the formula for $\psi_t(f)(x)$ on p.1380 
in \cite{Stannat03'}.) 
As a result we can show that 
\[ 
\limsup_{t\to\infty}\frac{1}{t}\log\Delta(t)
\ge -2\essinf b  
\] 
and accordingly 
$\ds{{\rm gap}(\overline{\cL_{1}}^{(2)}) \le \essinf b}$. 
In fact, the equality is valid 
because the opposite inequality is implied by  
(\ref{3.7}) with $\gamma=1$ 
with the help of Theorem 1.2 in \cite{Stannat05}. 
(See also Remark 2.15 in \cite{Stannat03}.) 
Thus ${\rm gap}(\overline{\cL_{\alpha}}^{(2)})$ 
is discontinuous at $\alpha=1$. \\ 
(ii) Theorem 4.3 would be regarded 
as a sort of continuous analogue of 
Theorem 2.6 in \cite{Liggett89}, 
which concerns a vector Markov process 
whose (countably many) components 
are independent Markov processes. 
\section{An application to generalized Fleming-Viot processes}
\setcounter{equation}{0} 
Throughout this section, we assume 
that $E$ is a compact metric space containing 
at least two distinct points. 
As mentioned in the Introduction, 
the previous results will be applied 
to study convergence to equilibrium 
for a class of generalized Fleming-Viot processes, 
whose state space is $\cM_1(E)$. 
These models have been discussed in \cite{H13}, 
where their stationary distributions were 
identified by exploiting connection with 
suitable measure-valued $\alpha$-CIR models. 
To be more precise, given $0<\alpha<1$ and 
$m\in\cM(E)^{\circ}$, we consider in this section 
the process associated with $\cL_{\alpha}$ 
in (\ref{2.1}) with $a\equiv 1\equiv b$ and 
the generalized Fleming-Viot process associated with 
\begin{eqnarray} 
\lefteqn{\cA_{\alpha}\Phi(\mu)}             \label{5.1} \\  
& := & 
\int_0^1\frac{B_{1-\alpha,1+\alpha}(du)}{u^2}\int_E\mu(dr)
\left[\Phi((1-u)\mu+u\delta_r)-\Phi(\mu)\right]  
                                             \nonumber \\ 
&  & 
+\int_0^1\frac{B_{1-\alpha,\alpha}(du)}{(\alpha+1)u}
\int_E m(dr)\left[\Phi((1-u)\mu+u\delta_r)-\Phi(\mu)\right], 
\qquad \mu\in\cM_1(E),                          \nonumber 
\end{eqnarray} 
where $B_{c_1,c_2}$ denotes the beta distribution 
with parameter $(c_1,c_2)$ and 
$\Phi$ belongs to the class $\cF_1$ of functions of the form 
$\Phi_f(\mu):=\lg\mu^{\otimes n},f \rg$ 
for some positive integer $n$ and $f\in C(E^n)$. 
It has been proved in \cite{H13} (Proposition 3.1) that 
the closure of $(\cA_{\alpha},\cF_1)$ generates 
a Feller semigroup $(S(t))_{t\ge0}$ on $C(\cM_1(E))$. 
Also, as observed in \cite{H13} (Proposition 3.3), 
the generators $\cL_{\alpha}$ and $\cA_{\alpha}$ together enjoy 
the following identity: 
\be 
\cL_{\alpha}\Psi(\eta)
=\Gamma(\alpha+2){\eta(E)^{-\alpha}}\cA_{\alpha}
\Phi\left(\eta(E)^{-1}\eta\right), 
\qquad \eta\in\cM(E)^{\circ}, 
                                               \label{5.2}  
\ee 
where $\Psi(\eta)=\Phi(\eta(E)^{-1}\eta)$ and 
$\Phi$ is in the linear span $\cF_0$ of functions of the form 
$\mu\mapsto\lg \mu, f_1\rg\cdots\lg \mu, f_n\rg$ 
with $f_i\in C_{++}(E)$ and $n$ being a positive integer. 
Notice that in one dimension 
(\ref{5.2}) takes the form (\ref{1.3}). 
Under the assumption that $m(E)>1$, 
it was proved in Proposition 3.4 of \cite{H13} that 
\be 
P_{\alpha}(\cdot)
:= 
\Gamma(\alpha+1)(m(E)-1)
E^{{Q}_{\alpha}}\left[\eta(E)^{-\alpha}; 
\eta(E)^{-1}\eta\in \cdot \right]
                                               \label{5.3}  
\ee 
is a unique stationary distribution of the process 
associated with $\cA_{\alpha}$, where $Q_{\alpha}$ 
is the stationary distribution of processes 
associated with (\ref{2.1}) with $a\equiv 1\equiv b$. 
The pre-factor $\Gamma(\alpha+1)(m(E)-1)$ 
on the right side arises as 
the normalizing constant. More generally, 
the following moment formula holds for the 
random variable $\eta(E)$ under $Q_{\alpha}$. 
\begin{lm}
Let $Q_{\alpha}$ be as in (\ref{2.17}) 
with $a\equiv 1 \equiv b$. Then 
it holds that 
\be 
E^{Q_{\alpha}}
\left[\eta(E)^{\beta}\right]
=
\left\{ 
\begin{array}{ll} 
\ds{\frac{\Gamma(1-\frac{\beta}{\alpha})
\Gamma(m(E)+\frac{\beta}{\alpha})}
{\Gamma(1-\beta)\Gamma(m(E))}}, 
& (-\alpha m(E)<\beta<\alpha) \\ 
\infty & (otherwise).  
\end{array} 
\right. 
                                               \label{5.4}  
\ee 
\end{lm}
{\it Proof.}~ 
Observe that, for each $u>0$, 
(\ref{2.17}) with $a\equiv 1 \equiv b$ and $f\equiv u$ reads  
\[ 
E^{{Q}_{\alpha}}\left[e^{-u\eta(E)}\right] 
=e^{-m(E)\log(1+u^{\alpha})} 
= (1+u^{\alpha})^{-m(E)}. 
\] 
We only need to consider the case $\beta<1$ 
because 
$E^{Q_{\alpha}}[\eta(E)^{\beta}]=\infty$ 
for each $\beta\ge 1$ is implied by 
$E^{Q_{\alpha}}[\eta(E)^{\alpha}]=\infty$. 
Since in this case 
$\Gamma(1-\beta)t^{-(1-\beta)}
=\int_0^{\infty}duu^{-\beta}e^{-u t}$ for $t>0$, 
we get with the help of Fubini's theorem  
\begin{eqnarray*} 
E^{{Q}_{\alpha}}\left[\eta(E)^{\beta}\right] 
& = & 
E^{{Q}_{\alpha}}\left[\eta(E)\eta(E)^{-(1-\beta)}\right]  \\ 
& = & 
\frac{1}{\Gamma(1-\beta)}
\int_0^{\infty}du u^{-\beta}
E^{{Q}_{\alpha}}\left[\eta(E)e^{-u\eta(E)}\right]  \\ 
& = & 
\frac{1}{\Gamma(1-\beta)}
\int_0^{\infty}du u^{-\beta} 
\left(-\frac{d}{du}
E^{{Q}_{\alpha}}\left[e^{-u\eta(E)}\right]\right)   \\ 
& = & 
\frac{\alpha m(E)}{\Gamma(1-\beta)}
\int_0^{\infty}du u^{-\beta} u^{\alpha-1} 
(1+u^{\alpha})^{-(m(E)+1)}                         \\ 
& = & 
\frac{m(E)}{\Gamma(1-\beta)}
\int_0^{\infty}dv v^{-\frac{\beta}{\alpha}}(1+v)^{-(m(E)+1)}, 
\end{eqnarray*} 
which is easily seen to diverge for 
$\beta\ge \alpha$ and $\beta\le -\alpha m(E)$. 
For $\beta\in (-\alpha m(E), \alpha)$, 
the above integral is equal to 
the right side of (\ref{5.4}). \qed 

\medskip

In the rest of the paper, we assume that $m(E)>1$. 
A key observation for the subsequent argument is 
that (\ref{5.2}) yields the following 
connection between Dirichlet forms: 
\be  
E^{{Q}_{\alpha}}
\left[(-\cL_{\alpha})\Psi\cdot\Psi\right]
=
\frac{\Gamma(\alpha+2)}{\Gamma(\alpha+1)(m(E)-1)}
E^{{P}_{\alpha}}\left[(-\cA_{\alpha})\Phi\cdot\Phi\right].  
                                                \label{5.5}  
\ee  
Here, the left side calls for some explanation 
since it can not necessarily be written as $\cE(\Psi)$. 
Instead, one can approximate such special $\Psi$'s 
suitably by functions in $D(\overline{\cL_{\alpha}}^{(2)})$. 
This rather technical point is 
the contents of the next lemma. 
\begin{lm}
Suppose that $m(E)>1$ and let 
$\cL_{\alpha}$ and $Q_{\alpha}$ be as above. 
Let $\Phi\in \cF_0$ be arbitrary and define 
$\Psi(\eta)=\Phi(\eta(E)^{-1}\eta)$ 
for $\eta\in\cM(E)^{\circ}$. 
Then there exists a sequence 
$\{\Psi_N\}\subset D(\overline{\cL_{\alpha}}^{(2)})$ 
such that $\Psi_N\to \Psi$ in $L^2(Q_{\alpha})$ 
and $\cE(\Psi_N)\to E^{{Q}_{\alpha}}
\left[(-\cL_{\alpha})\Psi\cdot\Psi\right]$ 
as $N\to\infty$. Moreover, 
${\rm var}_{{Q}_{\alpha}}(\Psi) 
\le 2E^{{Q}_{\alpha}}
\left[(-\cL_{\alpha})\Psi\cdot\Psi\right]$.  
\end{lm}
{\it Proof.}~ 
For simplicity we assume that  
$\Phi$ is of the form 
$\Phi(\mu)=\lg \mu, f_1\rg\cdots\lg \mu, f_n\rg$ 
for some $f_i\in C_{++}(E)$ and positive integer $n$. 
Accordingly  $\Psi(\eta)
=\lg \eta, f_1\rg\cdots\lg \eta, f_n\rg \eta(E)^{-n}$. 
Define $\{\Psi_N\}\subset \widetilde{\cF}$ by 
\[ 
\Psi_N(\eta)=\lg \eta, f_1\rg\cdots\lg \eta, f_n\rg 
(\eta(E)+1/N)^{-n}. 
\] 
In fact, $\Psi_N\in D(\overline{\cL_{\alpha}}^{(2)})$ and 
$\overline{\cL_{\alpha}}^{(2)}\Psi_N=\cL_{\alpha}\Psi_N$ 
as will be shown in the last half of the proof, 
and we temporarily suppose the validity of them. 
Obviously 
$\Psi_N\to \Psi$ boundedly and pointwise on $\cM(E)^{\circ}$ 
and so $\Psi_N\to \Psi$ in $L^2(Q_{\alpha})$. 
In view of Example preceding to Lemma 2.1 
and the calculations in the proof of it, 
one can verify that 
$\cL_{\alpha}\Psi_N(\eta) \to \cL_{\alpha}\Psi(\eta)$ 
for each $\eta\in\cM(E)^{\circ}$. 
Moreover, by virtue of Lemma 2.1 (ii) 
\[ 
|\cL_{\alpha}\Psi_N(\eta)|
\le C_1(\eta(E)+1/N)^{-\alpha}+C_2 
\] 
for some constants $C_1$ and $C_2$ 
independent of $\eta\in\cM(E)^{\circ}$ and $N$. 
Therefore, by Lebesgue's dominated convergence theorem 
\[ 
\cE(\Psi_N)
=
E^{{Q}_{\alpha}}\left[(-\overline{\cL_{\alpha}}^{(2)})
\Psi_N\cdot\Psi_N\right]
=
E^{{Q}_{\alpha}}\left[(-\cL_{\alpha})\Psi_N\cdot\Psi_N\right]
\to 
E^{{Q}_{\alpha}}\left[(-\cL_{\alpha})\Psi\cdot\Psi\right]
\] 
since $\eta(E)^{-\alpha}$ 
is integrable with respect to $Q_{\alpha}$ 
by Lemma 5.1 together with $m(E)>1$. 
In addition, 
${\rm var}_{{Q}_{\alpha}}(\Psi_N) 
\le 2\cE(\Psi_N)
=2E^{{Q}_{\alpha}}
\left[(-\cL_{\alpha})\Psi_N\cdot\Psi_N\right]$ 
by ${\rm gap}(\overline{\cL_{\alpha}}^{(2)})=1/2$. 
Letting $N\to\infty$ gives 
the required inequality for $\Psi$.

The rest of the proof is devoted to showing that 
$\Psi_N\in D(\overline{\cL_{\alpha}}^{(2)})$ 
and $\overline{\cL_{\alpha}}^{(2)}\Psi_N
=\cL_{\alpha}\Psi_N$ 
for arbitrarily fixed $N$. For this purpose, 
it suffices to construct 
$\{\Psi^{(k)}\}_{k=1}^{\infty}\subset \cF$ such that 
$\Psi^{(k)}\to \Psi_N$ and 
$\cL_{\alpha}\Psi^{(k)}\to \cL_{\alpha}\Psi_N$ 
boundedly and pointwise on $\cM(E)$ as $k\to\infty$. 
Clearly this reduces to finding a sequence 
$\{\varphi^{(k)}\}_{k=1}^{\infty}\subset C_0^2(\R^{n+1})$ 
which approximates to 
$\varphi(x_1,\ldots,x_n):=x_1\cdots x_n(x_{n+1}+1/N)^{-n}$ 
in some appropriate manner. 
To this end, take a sequence $\{\chi_k\}$ 
in $C_0^2(\R)$ with the following properties; 
$\{\chi_k'\}$ and $\{\chi_k''\}$ are uniformly bounded, 
\[  
0\le \chi_k(x)\le x \ \mbox{for all} \ x\in \R_+ 
\qquad 
\mbox{and} 
\qquad 
\chi_k(x)
=\left\{
\begin{array}{ll}
x, & (0\le x \le k) \\ 
0, & (x\ge 2k).  
\end{array} 
\right. 
\] 
Define 
$\widetilde{\varphi}^{(k)}\in C^2(\R_+^{n+1})$ by 
$\widetilde{\varphi}^{(k)}(x_1,\ldots,x_{n+1})
=\chi_k(x_1)\cdots\chi_k(x_n)(x_{n+1}+1/N)^{-n}$ 
and accordingly set 
\be 
\Psi^{(k)}(\eta)
=
\widetilde{\varphi}^{(k)} 
(\lg \eta,f_1\rg,\cdots,\lg \eta,f_n\rg,\eta(E)). 
                                                \label{5.6}  
\ee 
It follows that 
$\Vert \Psi^{(k)}\Vert_{\infty}\le\Vert \Psi_{N}\Vert_{\infty}$. 
While the support of 
$\widetilde{\varphi}^{(k)}$ is not compact, 
$\Psi^{(k)}(\eta)=0$ whenever 
$\eta(E)\ge 2k(\min_{1\le i\le n}\inf_{x\in E}f_i(x))^{-1}=:2kR$. 
Therefore, letting $\tau_k\in C_0^2(\R_+)$ 
be such that 
$\tau_k(x)=(x+1/N)^{-n}$ for any $x\in [0,2kR]$ 
and defining $\varphi^{(k)}\in C_0^2(\R_+^{n+1})$ by 
$\varphi^{(k)}(x_1,\ldots,x_{n+1})
=\chi_k(x_1)\cdots\chi_k(x_n)\tau_k(x_{n+1})$, 
we get $\Psi^{(k)}(\eta)
=\varphi^{(k)}(\lg \eta,f_1\rg,\cdots,\lg \eta,f_n\rg,\eta(E))$ 
for all $\eta\in \cM(E)$. 
Thus $\Psi^{(k)}\in \cF$ and 
$\Psi^{(k)}(\eta)=\Psi_N(\eta)$ 
whenever 
$\eta(E)<k(\max_{1\le i\le n}\sup_{x\in E} f_i(x))^{-1}(\le kR)$. 
Hence $\Psi^{(k)}\to \Psi_N$ boundedly and pointwise on $\cM(E)$ 
and $\cL_{\alpha}\Psi^{(k)}\to \cL_{\alpha}\Psi_N$ 
pointwise on $\cM(E)$. 
In addition, by Lemma 2.1 (ii) and 
analogous calculations for (\ref{5.6}) 
to Example preceding to Lemma 2.1 
one can show that 
$\{\cL_{\alpha}\Psi^{(k)}\}_{k=1}^{\infty}$ 
is uniformly bounded. 
This completes the proof of Lemma 5.2. \qed 

\medskip 

\noindent 
It follows from (\ref{5.5}) and Lemma 5.2 that 
for any $\Phi\in \cF_0$ 
\be 
{\rm var}_{{Q}_{\alpha}}(\Psi)
\le
\frac{2\Gamma(\alpha+2)}{\Gamma(\alpha+1)(m(E)-1)}
E^{{P}_{\alpha}}\left[(-\cA_{\alpha})\Phi\cdot\Phi\right], 
                                               \label{5.7}  
\ee 
where $\Psi(\eta)=\Phi(\eta(E)^{-1}\eta)$. 
Moreover, noting that by the 
Stone-Weierstrass theorem 
the linear span of functions $f$ on $E^n$ 
of the form 
$f(r_1,\ldots,r_n)=f_1(r_1)\cdots f_n(r_n)$ with 
$f_i\in C_+(E)$ is dense in $C(E^n)$, 
one can show, with the help of 
the expression (3.2) in \cite{H13} 
for $\cA_{\alpha}\Phi_f$ with $f\in C(E^n)$, 
that (\ref{5.7}) extends to any $\Phi\in \cF_1$. 
Indeed, that expression takes the form 
\[ 
\cA_{\alpha}\Phi_f(\mu) 
=
\lg\mu^{\otimes n},\Theta^{(n)}f\rg
+\lg\mu^{\otimes n},\Xi^{(n)}f\rg-c_n\Phi_f(\mu)  
\] 
for some nonnegative bounded operators 
$\Theta^{(n)}$, $\Xi^{(n)}:C(E^n)\to C(E^n)$ 
and some positive constant $c_n$ independent of 
$f$ and $\mu$, and so if 
$\{g_k\}\subset C(E^n)$, 
$g_k(r_1,\ldots,r_n)\to f(r_1,\ldots,r_n)$ 
uniformly on $E^n$ as $k\to\infty$, then 
$\cA_{\alpha}\Phi_{g_k}\to \cA_{\alpha}\Phi_f$ 
and $\Phi_{g_k}\to \Phi_f$ uniformly on $\cM_1(E)$ 
as $k\to\infty$.

In contrast, the variance functionals 
of $P_{\alpha}$ and of $Q_{\alpha}$, 
denoted by ${\rm var}_{P_{\alpha}}$ 
and ${\rm var}_{Q_{\alpha}}$ respectively, 
do not seem to enjoy any nice relation with 
each other. Although it is not clear if exponential 
convergence to equilibrium occurs for the 
process associated with $\cA_{\alpha}$, 
we are going to discuss a weaker ergodic property  
by introducing another functional 
${\rm osc}^2(\Phi)$ for $\Phi\in L^2(P_{\alpha})$, 
which is defined to be 
the essential supremum of the function 
\[ 
\cM_1(E)\times \cM_1(E)\ni (\mu,\mu') 
\mapsto |\Phi(\mu)-\Phi(\mu')|^2 
\]  
with respect to $P_{\alpha}\otimes P_{\alpha}$. 
Let $Z:\cM(E)^{\circ}\to \cM_1(E)$ be given by 
$Z(\eta)=\eta(E)^{-1}\eta$. Since by (\ref{5.3}) 
$P_{\alpha}$ and $Q_{\alpha}\circ Z^{-1}$ are 
mutually absolutely continuous, 
${\rm osc}^2(\Phi)$ coincides with 
the essential supremum of the function 
\[ 
\cM(E)^{\circ}\times \cM(E)^{\circ}\ni (\eta,\eta') 
\mapsto |(\Phi\circ Z)(\eta)-(\Phi\circ Z)(\eta')|^2 
\]  
with respect to $Q_{\alpha}\otimes Q_{\alpha}$. 
\begin{lm}
Suppose that $m(E)>1$ and 
let $Q_{\alpha}$ and $P_{\alpha}$ be as in (\ref{5.3}). 
Let $\Phi\in L^2(P_{\alpha})$ 
be arbitrary and define 
$\Psi(\eta)=\Phi(\eta(E)^{-1}\eta)$ 
for $\eta\in\cM(E)^{\circ}$. Then for any $q>1$ 
\be 
{\rm var}_{P_{\alpha}}(\Phi) 
\le 
\Gamma(\alpha+1)(m(E)-1)\left(E^{Q_{\alpha}}
\left[\eta(E)^{-\alpha q}\right]
{\rm osc}^2(\Phi)\right)^{1/q}
\left({\rm var}_{Q_{\alpha}}(\Psi) \right)^{1/p}, 
                                               \label{5.8}
\ee 
where $p>1$ is such that $1/p+1/q=1$. 
\end{lm}
{\it Proof.}~
It follows from (\ref{5.3}) that 
\begin{eqnarray*} 
{\rm var}_{P_{\alpha}}(\Phi) 
& = & 
\inf_{c\in \R}
E^{P_{\alpha}}
\left[\left(\Phi-c\right)^2\right]  
\ \le \ 
E^{P_{\alpha}}
\left[\left(\Phi-E^{Q_{\alpha}}[\Psi]\right)^2\right]   \\ 
& = & 
\Gamma(\alpha+1)(m(E)-1)
E^{{Q}_{\alpha}}\left[\eta(E)^{-\alpha}
\left(\Psi(\eta)-E^{Q_{\alpha}}[\Psi]\right)^2\right].  
\end{eqnarray*} 
By H\"older's inequality 
\begin{eqnarray*} 
\lefteqn{
E^{{Q}_{\alpha}}\left[\eta(E)^{-\alpha}
\left(\Psi(\eta)-E^{Q_{\alpha}}[\Psi]\right)^2\right]} \\  
& \le & 
\left(E^{{Q}_{\alpha}}\left[\eta(E)^{-\alpha q}
\left|\Psi(\eta)-E^{Q_{\alpha}}[\Psi]\right|^{2}
\right]\right)^{1/q}
\left(E^{{Q}_{\alpha}}\left[
\left|\Psi(\eta)-E^{Q_{\alpha}}[\Psi]\right|^{2}
\right]\right)^{1/p}                                \\ 
& \le  & 
\left(E^{{Q}_{\alpha}}\left[\eta(E)^{-\alpha q}
\right]{\rm osc}^2(\Phi)\right)^{1/q}
\left({\rm var}_{Q_{\alpha}}(\Psi)\right)^{1/p},  
\end{eqnarray*} 
where the last inequality is seen from 
\begin{eqnarray*} 
\lefteqn{ 
E^{{Q}_{\alpha}}\left[\eta(E)^{-\alpha q}
\left|\Psi(\eta)-E^{Q_{\alpha}}[\Psi]\right|^{2}\right]} \\  
& \le & 
E^{{Q}_{\alpha}\otimes {Q}_{\alpha}}
\left[\eta(E)^{-\alpha q}
\left|\Psi(\eta)-\Psi(\eta')\right|^{2}\right]   \\ 
& = & 
E^{{Q}_{\alpha}\otimes {Q}_{\alpha}}
\left[\eta(E)^{-\alpha q}
\left|(\Phi\circ Z)(\eta)-(\Phi\circ Z)(\eta')
\right|^{2}\right]. 
\end{eqnarray*}
These inequalities together prove (\ref{5.8}). \qed 

\medskip 

Our strategy is to carry out the well-known procedure 
to show algebraic convergence to equilibrium. 
More specifically, our ingredient for this 
is Theorem 2.2 in \cite{Liggett91}. 
(See also \cite{RW}.) 
Let $\{S^2(t)\}_{t\ge 0}$ be the strongly continuous 
semigroup on $L^2(P_{\alpha})$ 
of the process associated with $\cA_{\alpha}$. 
To be more precise, $\{S^2(t)\}_{t\ge 0}$ is 
defined to be the $C_0$-semigroup on $L^2(P_{\alpha})$ 
generated by the closure 
$(\overline{\cA_{\alpha}}^{(2)},
D(\overline{\cA_{\alpha}}^{(2)}))$ 
of $(\cA_{\alpha}, \cF_1)$ as 
an operator on $L^2(P_{\alpha})$. 
As in the case of $({T}^2(t))_{t\ge0}$, 
$({S}^2(t))_{t\ge0}$ on $L^2(P_{\alpha})$ 
coincides with $(S(t))_{t\ge0}$ 
when restricted to $C(\cM_1(E))$. 
The following property of $({S}^2(t))_{t\ge0}$ 
is needed for the abovementioned strategy. 
\begin{lm}
Let $({S}^2(t))_{t\ge0}$ be as above. 
Then it holds that 
${\rm osc}^2(S^2(t)\Phi)\le {\rm osc}^2(\Phi)$ 
for all $t>0$ and $\Phi\in L^2(P_{\alpha})$. 
\end{lm}
{\it Proof.}~ 
We may and do assume additionally 
that ${\rm osc}^2(\Phi)<\infty$. 
Let $\cS$ be an arbitrary 
Borel set of $\cM_1(E)\times \cM_1(E)$ and 
$P_t(\mu,\cdot)$ $(t>0, \mu\in\cM_1(E))$ 
denote the transition function of the process 
associated with $\cA_{\alpha}$. 
It then follows that for any $q>1$ 
\begin{eqnarray*} 
\lefteqn{ 
E^{{P}_{\alpha}\otimes {P}_{\alpha}} 
\left[
\left|S^2(t)\Phi(\mu)-S^2(t)\Phi(\mu')\right|^2
\one_{\cS}(\mu,\mu')\right]}                       \\ 
& \le & 
E^{{P}_{\alpha}\otimes {P}_{\alpha}} 
\left[
\int_{\cM_1(E)}P_t(\mu,d\nu)\int_{\cM_1(E)}P_t(\mu',d\nu')
\left|\Phi(\nu)-\Phi(\nu')\right|^2
\one_{\cS}(\mu,\mu')\right]                        \\ 
& \le & 
\left(E^{{P}_{\alpha}\otimes {P}_{\alpha}} 
\left[
\int_{\cM_1(E)}P_t(\mu,d\nu)\int_{\cM_1(E)}P_t(\mu',d\nu')
\left|\Phi(\nu)-\Phi(\nu')\right|^{2q}\right]\right)^{1/q} \\ 
& & 
\times \left(({P}_{\alpha}\otimes {P}_{\alpha})(\cS)\right)^{1/p} \\ 
& = & 
\left(E^{{P}_{\alpha}\otimes {P}_{\alpha}} 
\left[\left|\Phi(\mu)-\Phi(\mu')\right|^{2q}\right]\right)^{1/q} 
\left(({P}_{\alpha}\otimes {P}_{\alpha})(\cS)\right)^{1/p} \\ 
& \le & 
{\rm osc}^2(\Phi)
\left(({P}_{\alpha}\otimes {P}_{\alpha})(\cS)\right)^{1/p}, 
\end{eqnarray*} 
where $p>1$ is such that $1/p+1/q=1$ and 
the equality is implied by the stationarity of $P_{\alpha}$. 
By letting $q\to \infty$ or $p\downarrow 1$ 
\[ 
E^{{P}_{\alpha}\otimes {P}_{\alpha}} 
\left[\left|S^2(t)\Phi(\mu)-S^2(t)\Phi(\mu')\right|^{2}
\one_{\cS}(\mu,\mu')\right] 
\le 
{\rm osc}^2(\Phi)({P}_{\alpha}\otimes {P}_{\alpha})(\cS). 
\] 
Since $\cS$ is arbitrary, we conclude that 
${\rm osc}^2(S^2(t)\Phi)\le {\rm osc}^2(\Phi)$.     \qed 

\medskip 

\noindent 
We can at last state the main result of this section. 
\begin{th}
Suppose that $m(E)>1$. Let $P_{\alpha}$ be as in (\ref{5.3}) 
and $({S}^2(t))_{t\ge0}$ be as above. 
Then for any $\Phi\in L^2(P_{\alpha})$ 
\begin{eqnarray} 
\lefteqn{\limsup_{t\to \infty}
{\rm var}_{P_{\alpha}}(S^2(t)\Phi)
\frac{t^{m(E)-1}}{\log t} }                  \nonumber \\ 
& \le & 
\frac{e\Gamma(\alpha)(m(E)-1)}{\Gamma(\alpha m(E))}
\left(\Gamma(\alpha+2)(m(E)-1)\right)^{m(E)-1}
{\rm osc}^2(\Phi).  
                                               \label{5.9}  
\end{eqnarray} 
\end{th}
{\it Proof.}~Given $\Phi\in L^2(Q_{\alpha})$, 
let $\Psi(\eta)=\Phi(\eta(E)^{-1}\eta)$ 
for $\eta\in\cM(E)^{\circ}$. 
We claim that (\ref{5.7}) is extended as 
\be 
{\rm var}_{{Q}_{\alpha}}(\Psi)
\le
\frac{2\Gamma(\alpha+2)}{\Gamma(\alpha+1)(m(E)-1)}
E^{{P}_{\alpha}}\left[(-\overline{\cA_{\alpha}}^{(2)})
\Phi\cdot\Phi\right] 
                                               \label{5.10}  
\ee 
for any $\Phi\in D(\overline{\cA_{\alpha}}^{(2)})$. 
Indeed, there exists a sequence $\{\Phi_N\}\subset \cF_1$ 
such that $\Phi_N\to \Phi$ and 
$\cA_{\alpha}\Phi_N\to \overline{\cA_{\alpha}}^{(2)}\Phi$ 
in $L^2(P_{\alpha})$ as $N\to\infty$, 
and according to (\ref{5.7}) for elements of $\cF_1$ 
we have for each $N=1,2,\ldots$ 
\be 
E^{Q_{\alpha}}[(\Psi_N)^2]
-\left(E^{Q_{\alpha}}[\Psi_N]\right)^2  
\le \frac{2\Gamma(\alpha+2)}{\Gamma(\alpha+1)(m(E)-1)}
E^{{P}_{\alpha}}
\left[(-\cA_{\alpha})\Phi_N\cdot\Phi_N\right], 
                                               \label{5.11}  
\ee 
where $\Psi_N(\eta)=\Phi_N(\eta(E)^{-1}\eta)$. 
Taking a subsequence if necessary, 
we can assume that $\Phi_N\to \Phi$ $P_{\alpha}$-a.s. 
and so $\Psi_N\to \Psi$ $Q_{\alpha}$-a.s. 
Since (\ref{5.11}) implies that 
$\{\Psi_N\}$ is bounded in $L^2(Q_{\alpha})$, 
$\Psi\in L^2(Q_{\alpha})$ and 
$\Psi_N\to \Psi$ in $L^1(Q_{\alpha})$. 
Letting $N\to\infty$ in (\ref{5.11}) yields 
\[ 
E^{Q_{\alpha}}[\Psi^2]
-\left(E^{Q_{\alpha}}[\Psi]\right)^2  
\le \frac{2\Gamma(\alpha+2)}{\Gamma(\alpha+1)(m(E)-1)}
E^{{P}_{\alpha}}
\left[(-\overline{\cA_{\alpha}}^{(2)})
\Phi\cdot\Phi\right] 
\] 
with the help of Fatou's lemma. 
Thus (\ref{5.10}) holds for 
$\Phi\in D(\overline{\cA_{\alpha}}^{(2)})$.

Let $q\in(1,m(E))$ be arbitrary. 
Combining (\ref{5.10}) with (\ref{5.8}) leads to  
\[ 
{\rm var}_{P_{\alpha}}(\Phi) 
\le 
\left(\cV(\Phi)\right)^{1/q}
\left(2\Gamma(\alpha+2)
E^{{P}_{\alpha}}\left[(-\overline{\cA_{\alpha}}^{(2)})
\Phi\cdot\Phi\right]\right)^{1/p}, 
\quad \Phi\in D(\overline{\cA_{\alpha}}^{(2)}),  
\]  
where 
$\cV(\Phi)=\Gamma(\alpha+1)(m(E)-1)E^{Q_{\alpha}}
\left[\eta(E)^{-\alpha q}\right]{\rm osc}^2(\Phi)$. 
Thus the condition (2.3) of Theorem 2.2 in 
\cite{Liggett91} is fulfilled with 
a quadratic functional $\Phi\mapsto \cV(\Phi)$ 
and $C:=(2\Gamma(\alpha+2))^{1/p}$. 
In addition, by Lemma 5.4 
$\cV(S^2(t)\Phi)\le \cV(\Phi)$ 
for all $t\ge 0$ and $\Phi\in L^2(P_{\alpha})$. 
Therefore, by the assertion (i) 
of Theorem 2.2 in \cite{Liggett91} together 
with the calculation in the final part of its proof, 
we obtain  
\[ 
{\rm var}_{P_{\alpha}}(S^2(t)\Phi)
\le  
\cV(\Phi)C^q\left(\frac{q/p}{2t}\right)^{q/p} 
\] 
for all $\Phi\in L^2(P_{\alpha})$ and $t>0$. 
Because of $q/p=q-1$, this is rewritten as 
\begin{eqnarray} 
\lefteqn{{\rm var}_{P_{\alpha}}(S^2(t)\Phi)}  \nonumber \\ 
& \le & 
\Gamma(\alpha+1)(m(E)-1)
E^{Q_{\alpha}}\left[\eta(E)^{-\alpha q}\right]
{\rm osc}^2(\Phi) 
\left(\frac{\Gamma(\alpha+2)(q-1)}{t}\right)^{q-1}.  
                                               \label{5.12}  
\end{eqnarray}

As the final step, we will optimize the value 
of $q$ in (\ref{5.12}) for each $t>0$ large enough. 
Observe from (\ref{5.4}) 
that the right side of (\ref{5.12}) becomes  
\be 
\frac{\Gamma(\alpha)\Gamma(q)\Gamma(m(E)-q)}
{\Gamma(\alpha q)\Gamma(m(E)-1)}
\left(\frac{\Gamma(\alpha+2)(q-1)}{t}\right)^{q-1}
{\rm osc}^2(\Phi) .    
                                               \label{5.13}  
\ee 
For $t$ sufficiently large, take 
$q=q(t):=m(E)-\delta(t)\in (1,m(E))$, 
where $\delta(t):=1/(\log t)$ is verified to 
minimize the function $0<\delta\mapsto t^{\delta}/\delta$. 
Then, noting that 
\[ 
\Gamma(m(E)-q(t))
=\frac{\Gamma(\delta(t)+1)}{\delta(t)}
=\Gamma(\delta(t)+1)\cdot \log t 
\] 
and $t^{\delta(t)}=e$, we see that (\ref{5.13}) 
with this choice of $q$ equals 
\[ 
\frac{\Gamma(\alpha)\Gamma(q(t))\Gamma(\delta(t)+1)}
{\Gamma(\alpha q(t))\Gamma(m(E)-1)}\cdot
\frac{\log t}{t^{m(E)-1}}\cdot e 
\left(\Gamma(\alpha+2)(q(t)-1)\right)^{q(t)-1}
{\rm osc}^2(\Phi).  
\] 
This upper bound for ${\rm var}_{P_{\alpha}}(S^2(t)\Phi)$ 
immediately gives (\ref{5.9}).  
The proof of Theorem 5.5 is complete. \qed 

\medskip 

\noindent 
What is unpleasant to us is that 
we do not know whether 
${\rm gap}(\overline{\cA_{\alpha}}^{(2)})=0$ or not. 
One difficulty is that 
any useful expression for 
the variance functional with respect to $P_{\alpha}$ 
nor Dirichlet form associated with $\cA_{\alpha}$ 
does not seem available at least for 
conventional choice of test functions. 
Besides, our argument in this section 
does not work in the case where $0<m(E)\le 1$ 
although the process associated with $\cA_{\alpha}$ 
still has a unique stationary distribution. 
(See Theorem 3.2 in \cite{H13}.) 

\bigskip 

\noindent 
{\bf Acknowledgment.}~ 
The author is grateful to an anonymous 
referee for helpful comments that improved 
the presentation of the paper. 
Part of the work was carried out under 
the ISM Cooperative Research Program 
(2013$\cdot$ISM$\cdot$CRP-5010).

\end{document}